\documentclass[english]{IEEEtran}
\usepackage[T1]{fontenc}
\usepackage[latin9]{inputenc}
\usepackage{pifont}
\usepackage{amsthm}
\usepackage{amsmath}
\usepackage{amssymb}
\usepackage{graphicx}
\usepackage{epstopdf}
\usepackage{paralist}
\usepackage[T1]{fontenc}
\usepackage{epsfig} 
\usepackage{cite}
\usepackage{dsfont}
\usepackage{bbm}
\usepackage[small]{caption}
\usepackage{subfigure}	
\usepackage{float}
\usepackage{graphicx}

\makeatletter

\theoremstyle{plain} \newtheorem{thm}{\protect\theoremname}
\theoremstyle{plain} \newtheorem{lem}{\protect\lemmaname}
\theoremstyle{plain} \newtheorem{cor}{\protect\corollaryname}
\theoremstyle{plain} 
\theoremstyle{definition} 

\providecommand{\lemmaname}{Lemma}
\providecommand{\theoremname}{Theorem}
\providecommand{\examplename}{Example}
\providecommand{\corollaryname}{Corollary}

\newcommand{\subscr}[2]{#1_{\textup{#2}}}										

\@ifundefined{showcaptionsetup}{}{
 \PassOptionsToPackage{caption=false}{subfig}}
\usepackage{subfig}
\AtBeginDocument{
  
}

\makeatother

\usepackage{babel}

\makeatletter
\adddialect\l@ENGLISH\l@english
\makeatother

\begin{document}
\title{Synchronization of Nonlinear Circuits in Dynamic Electrical Networks with General Topologies}

\author{Sairaj V. Dhople, \emph{Member, IEEE}, Brian B. Johnson, \emph{Member, IEEE}, \linebreak Florian D\"{o}rfler, \emph{Member, IEEE}, and Abdullah Hamadeh
\thanks{S. V. Dhople is with the Department of Electrical and Computer Engineering at the University of Minnesota, Minneapolis, MN (email: \texttt{sdhople@UMN.EDU}); B. B. Johnson is with the Power Systems Engineering Center at the National Renewable Energy Laboratory, Golden, CO (email: \texttt{brian.johnson@NREL.GOV}); F. D\"{o}rfler is with the Department of Electrical Engineering at the University of California at Los Angeles, Los Angeles, CA (email: \texttt{dorfler@SEAS.UCLA.EDU}); A. O. Hamadeh is with the Department of Mechanical Engineering, Massachusetts Institute of Technology, Cambridge, MA (e-mail: \texttt{ahamadeh@MIT.EDU}).}}

\maketitle

\begin{abstract}
Sufficient conditions are derived for global asymptotic synchronization in a system of identical nonlinear electrical circuits coupled through linear time-invariant (LTI) electrical networks. In particular, the conditions we derive apply to settings where: i) the nonlinear circuits are composed of a parallel combination of passive LTI circuit elements and a nonlinear voltage-dependent current source with finite gain; and ii) a collection of these circuits are coupled through either uniform or homogeneous LTI electrical networks. Uniform electrical networks have identical per-unit-length impedances. Homogeneous electrical networks are characterized by having the same effective impedance between any two terminals with the others open circuited. Synchronization in these networks is guaranteed by ensuring the stability of an equivalent coordinate-transformed differential system that emphasizes signal differences. The applicability of the synchronization conditions to this broad class of networks follows from leveraging recent results on structural and spectral properties of Kron reduction---a model-reduction procedure that isolates the interactions of the nonlinear circuits in the network. The validity of the analytical results is demonstrated with simulations in networks of coupled Chua's circuits.  
\end{abstract}

\section{Introduction}

\IEEEPARstart{S}{ynchronization} of nonlinear electrical circuits coupled through complex networks is integral to modeling, analysis, and control in application areas such as the ac electrical grid, solid-state circuit oscillators, semiconductor laser arrays, secure communications, and microwave oscillator arrays~\cite{Strogatz_Book01,Dorfler-13-Synch}. This paper focuses on the global asymptotic synchronization of terminal voltages in a class of nonlinear circuits coupled through passive LTI electrical networks. We assume that the nonlinear circuits are composed of a parallel combination of a passive LTI circuit and a nonlinear voltage-dependent current source with finite gain. A collection of such identical circuits are coupled through uniform or homogeneous passive LTI electrical networks. Uniform networks have identical per-unit-length impedances and include purely resistive and lossless networks as special cases. Homogeneous electrical networks are characterized by identical effective impedances between the terminals (essentially, the impedance between any two terminals with the others open circuited). Section~\ref{sec:GASSufficient1} provides precise definitions of these network types.

The nonlinear-circuit models, and the uniform and homogeneous networks examined in this work offer a broad level of generality and ensure a wide applicability of the analytical results to many settings. For instance, a variety of chaotic and hyperchaotic circuits as well as nonlinear oscillators~\cite{Chua-1990-Canonical,Gismero-1990,Chua_Passivity99,Kyprianidis-1999,Malhame-2000,Kyprianidis-2001,Torres_Power12,Torres_Power13,Johnson_TCAS13} can be modeled as a parallel connection of a linear subsystem and a nonlinear voltage-dependent current source with finite gain. Similarly, the types of networks that our results accommodate, facilitate the analysis of varied interconnections between the nonlinear circuits. In general, we study interconnecting networks that are \emph{dynamic}, i.e., the network can contain capacitive or inductive storage elements. For uniform networks, the per-unit-length line impedances can be complex (i.e., not exclusively resistive or reactive) and the network topology can be arbitrary. With regard to resistive (lossless) networks, we allow the branch resistances (respectively, reactances) and the network topology to be arbitrary.  Finally, homogeneous networks are frequently encountered in symmetric engineered setups (e.g., power grid monitoring and electrical impedance tomography), in large random networks or regular lattices, as well as in idealized settings where all terminals are electrically uniformly distributed with respect to each other~\cite{Dorfler-13}.

The analytical approach adopted in this paper builds on previous work in~\cite{Hamadeh_Constructive10,Hamadeh_Thesis,Hamadeh_Designing12,Johnson_TCAS13}, where $\mathcal{L}_{2}$ methods were used to  analyze synchronization in feedback systems thereby offering an alternate perspective compared to a rich body of literature that has examined synchronization problems with Lyapunov- and passivity-based mehods~\cite{Wu-1995,Wu-1998, Pogromsky_Cooperative01,GBS_Thesis,Wu-2005,Arcak_Passivity07, GBS_Analysis07, Hamadeh_Global12,Torres_Power12,Torres_Power13}. To investigate synchronization, the linear and nonlinear subsystems in the network of coupled nonlinear electrical circuits are compartmentalized, and a coordinate transformation is applied to recover an equivalent \emph{differential system} that emphasizes signal differences. Once the differential $\mathcal{L}_2$ gains of the linear and nonlinear subsystems are identified, synchronization can be guaranteed by ensuring the stability of the coordinate-transformed differential system with a small-gain argument. 

The suite of synchronization conditions presented in this paper generalize our previous efforts in~\cite{Johnson_TCAS13} (which were limited to electrical networks with a star topology) to arbitrary network topologies. Integral to the analysis that allows us to investigate varied topologies is a model-reduction procedure called \emph{Kron reduction}~\cite{Kron-1939}. This procedure explicitly uncovers the interactions between the nonlinear electrical circuits, while systematically eliminating exogenous nodes in the network. A major contribution of this work pertains to leveraging recent results on structural and spectral properties of Kron reduction~\cite{Dorfler-13} in deriving synchronization conditions. Towards this end, another significant contribution is that some key lemmas from~\cite{Dorfler-13} are extended from the real-valued and symmetric to the complex-symmetric (and not necessarily Hermitian) case, and we also offer converse results to some statements in~\cite{Dorfler-13}.

The remainder of this manuscript is organized as follows. Section~\ref{sec:Preliminaries} establishes some mathematical preliminaries and notation. In Section~\ref{sec:SystemDescription}, we describe the nonlinear electrical circuits, and describe their network interactions by constructing the electrical admittance matrix that couples them. In Section~\ref{sec:Statement}, we formulate the problem statement, and derive the differential system. Synchronization conditions for networks with and without shunt elements are then derived in Sections~\ref{sec:GASSufficient1} and~\ref{sec:GASSufficient2}, respectively. Simulations-based case studies are provided in Section~\ref{sec:CaseStudies} to validate the approach. We conclude the paper in Section~\ref{sec:Conclusions} by highlighting a few pertinent directions for future work. 

\section{Notation and Preliminaries\label{sec:Preliminaries}}

Given a complex-valued $N$-tuple $\{u_1, \ldots, u_N\}$, denote the corresponding column vector as $u = [u_1, \ldots, u_N]^\mathrm{T}$, where 
$(\cdot)^\mathrm{T}$ denotes transposition (without conjugation). Denote the $N \times N$ identity matrix as $I$, and the $N$-dimensional vectors of all ones and zeros as $\mathbf{1}$ and $\mathbf{0}$, respectively. The Moore-Penrose pseudo inverse of a matrix $U$ is denoted by $U^\dagger$. Let $\mathrm{j} = \sqrt{-1}$ be the imaginary unit. Denote the Laplace transform of a continuous-time function $f(t)$ by $f$. Cardinality of the set $\mathcal{N}$ is denoted by $|\mathcal{N}|$.

The Euclidean norm of a complex vector, $u$, is denoted by $\left\Vert u\right\Vert _2$ and is defined as 
\begin{equation}
\left\Vert u\right\Vert _2:=\sqrt{u^{*}u},\label{eq:EuclidianNormDef}
\end{equation}
where $(\cdot)^{*}$ signifies the conjugate transpose. The space of all piecewise continuous functions such that
\begin{equation}
\left\Vert u\right\Vert _{\mathcal{L}_{2}}:=\sqrt{\intop_{0}^{\infty}u\left(t\right)^{\mathrm{T}}u\left(t\right)\,\mathrm{d}t} < \infty,
\end{equation}
is denoted as $\mathcal{L}_2$, where $\left\Vert u\right\Vert_{\mathcal{L}_2}$ is referred to as the \emph{$\mathcal{L}_{2}$ norm} of $u$ \cite{VanDerSchaft_Book96}. If $u \in \mathcal{L}_2$, then $u$ is said to be \emph{bounded}.

A causal system, $\mathcal{H}$, with input $u$ and output $y$ is \emph{finite-gain $\mathcal{L}_{2}$ stable} if there exist finite and non-negative constants $\gamma$ and $\eta$ such that
\begin{equation}
  \left\Vert y\right\Vert _{\mathcal{L}_{2}}=:\left\Vert \mathcal{H}\left(u\right)\right\Vert _{\mathcal{L}_{2}}\leq\gamma\left\Vert u\right\Vert _{\mathcal{L}_{2}}+\eta,\quad\forall u\in\mathcal{L}_{2}.\label{eq:Finite_gain_stableDef}
\end{equation}
The smallest value of $\gamma$ for which there exists an $\eta$ such that~\eqref{eq:Finite_gain_stableDef} is satisfied is called the \emph{$\mathcal{L}_{2}$} \emph{gain} of the system. If $\mathcal{H}$ is linear and can be represented by the transfer matrix $H: \mathbb C \to \mathbb{C}^{N\times N}$, it can be shown that the $\mathcal{L}_{2}$ gain of $\mathcal{H}$ is equal to its H-\emph{infinity norm}, denoted by $\left\Vert \mathcal{H}\right\Vert _{\infty}$, and defined as
\begin{equation}
\gamma\left(\mathcal{H}\right)=\left\Vert \mathcal{H}\right\Vert _{\infty}:=\underset{\omega\in\mathbb{R}}{\mathrm{sup}}\,\frac{\left\Vert H\left(\mathrm{j}\omega\right)u\left(\mathrm{j}\omega\right)\right\Vert _{2}}{\left\Vert u\left(\mathrm{j}\omega\right)\right\Vert _{2}},\label{eq:InfNormSingularValue}
\end{equation}
where $\left\Vert u\left(\mathrm{j}\omega\right)\right\Vert _{2}=1$, provided that all poles of $H$ have strictly negative real parts~\cite{Khalil_Book02}. For a single-input single-output transfer function $h: \mathbb C \to \mathbb{C}$, $\gamma\left(\mathcal{H}\right)=\left\Vert
  \mathcal{H}\right\Vert
_{\infty}=\underset{\omega\in\mathbb{R}}{\mathrm{sup}}\,\left\Vert
  h\left(\mathrm{j}\omega\right)\right\Vert _{2}$.

A construct we will find particularly useful in assessing signal differences is the $N\times N$ \emph{projector matrix}~\cite{GBS_Analysis07,GBS_Thesis,Hamadeh_Thesis}, which is denoted by $\Pi$, and defined as 
\begin{equation}
  \label{eq:ProjectorDefinition}
  \Pi := I - \frac{1}{N}\mathbf{11}^\mathrm{T}.
\end{equation}
For a vector $u$, we define $\widetilde{u} := \Pi u$ to be the corresponding \emph{differential vector} \cite{GBS_Analysis07,GBS_Thesis,Hamadeh_Constructive10, Hamadeh_Thesis, Hamadeh_Designing12}. 

Given a symmetric and nonnegative matrix $A \in \mathbb R^{N \times N}$ associated with an undirected and weighted graph, we define its {\em Laplacian matrix} $L$ component-wise by $l_{nm} = -a_{nm}$ for off-diagonal elements and $l_{nn}  = \sum_{m=1}^{N} a_{nm}$ for diagonal elements. The Laplacian matrix has zero row and column sums, it is symmetric and positive semidefinite, and its zero eigenvalue is simple if and only if the graph is connected. The Laplacian of a complete graph with unit weights is 
\begin{equation} 
\Gamma = N I - \mathbf{1}\mathbf{1}^{\mathrm{T}} = N \Pi.
\end{equation}

A causal system with input $u$ and output $y$ is said to be \emph{differentially finite $\mathcal{L}_2$ gain stable} if there exist finite, non-negative constants, $\widetilde{\gamma}$ and $\widetilde{\eta}$, such that
\begin{equation}
  \label{eq:DiffL2GainStability}
  \left\Vert \widetilde{y}\right\Vert _{\mathcal{L}_{2}}\leq\widetilde{\gamma}\left\Vert \widetilde{u}\right\Vert _{\mathcal{L}_{2}}+\widetilde{\eta},\quad\forall\,\widetilde{u}\in\mathcal{L}_{2},
\end{equation}
where $\widetilde{y} = \Pi y$. The smallest value of $\widetilde{\gamma}$ for which there exists a non-negative value of $\widetilde{\eta}$ such that~\eqref{eq:DiffL2GainStability} is satisfied is referred to as the \emph{differential $\mathcal{L}_2$ gain} of $H$. The differential $\mathcal{L}_2$ gain of a system provides a measure of the largest amplification imparted to input signal differences.

Consider two systems that are modeled by transfer matrices $A$ and $B$. The \emph{linear fractional transformation} is the transfer matrix of the negative-feedback interconnection of these systems, and it is given by~\cite{Zhou96robust_book}
\begin{equation}
\mathcal{F}\left(A(s),B(s) \right) := \left(I +  A(s)B(s) \right)^{-1}A(s).
\label{eq:LFT}
\end{equation}

For an electrical network with admittance matrix $Y$, the {\em effective impedance} $z_{nm}$ between nodes $n$ and $m$ is the potential difference between nodes $n$ and $m$, when a unit current is injected in node $n$ and extracted from node $m$. 
In this case, the current-balance equations are $e_{n} - e_{m} = Y \upsilon$, where $e_n$ is the canonical\,vector of all zeros except with a 1 in the $n^\mathrm{th}$ position, and $\upsilon$ is the vector of the resulting nodal voltages. The effective impedance is then
\begin{equation}
z_{nm} = (e_{n} - e_{m})^\mathrm{T} \upsilon = (e_{n} - e_{m})^\mathrm{T} Y^{\dagger} (e_n - e_m) \,. \label{eq:zeff}
\end{equation}
The effective impedance is an electric and graph-theoretic distance measure, see~\cite{Dorfler-13} for details and further references.

\section{System of Coupled Nonlinear Electrical Circuits}\label{sec:SystemDescription}
We begin this section with a brief description of the type of nonlinear electrical circuits for which we derive sufficient synchronization conditions. Next, we describe the electrical network that couples the nonlinear electrical circuits. 

\subsection{Nonlinear Circuit Model\label{sub:System-Description}}
\begin{figure}[b]
\begin{centering}
\includegraphics{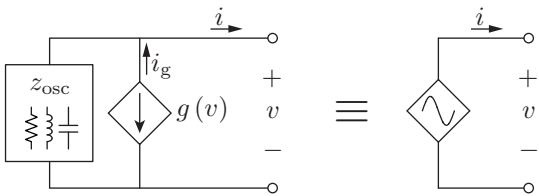}
\par\end{centering}
\caption{Electrical schematic of the nonlinear circuit studied in this work. Each circuit is composed of a linear subsystem modeled by a passive impedance, $z_\mathrm{osc}$, and a nonlinear voltage-dependent current source, $g(\cdot)$. Circuit symbol used to represent the nonlinear circuit is shown on the right. \label{Fig:SingleOscillator}}
\end{figure}

\begin{figure*}[t]
\begin{center}
\subfigure[]{\includegraphics[width=0.28\linewidth]{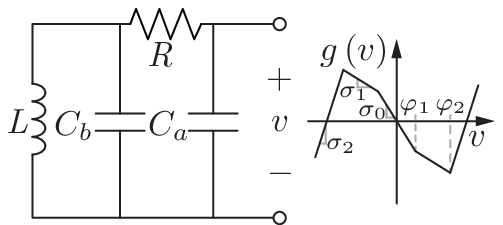}}\quad \quad \quad
\subfigure[]{\includegraphics[width=0.27\linewidth]{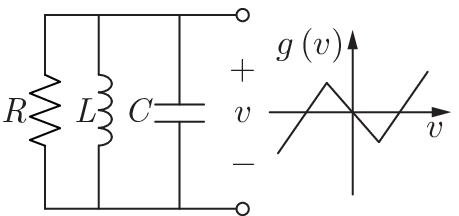}}\quad \quad \quad
\subfigure[]{\includegraphics[width=0.28\linewidth]{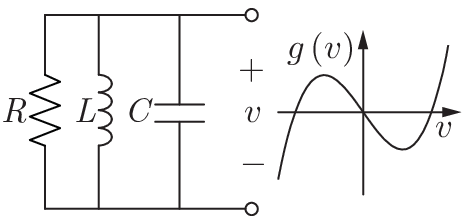}}
\end{center}
\caption{The linear-subsystem impedance, $z_\mathrm{osc}$, and the nonlinear voltage-dependent current source, $g(\cdot)$, illustrated for (a) Chua's circuit, (b) the Dead-zone oscillator, and (c) the Van der Pol oscillator.}
\label{Fig:Examples}
\end{figure*}

An electrical schematic of the nonlinear circuits studied in this work is depicted in Fig.~\ref{Fig:SingleOscillator}. Each circuit has a linear subsystem composed of an arbitrary connection of passive circuit elements described by the impedance, $z_\mathrm{osc}\in \mathbb{C}$, and a nonlinear voltage-dependent current source $i_\mathrm{g}=-g(v)$. We will require that the maximum slope of the function $g(\cdot)$ be bounded:
\begin{equation}
\sigma := \sup_{v\in\mathbb{R}}\left|\frac{\text{d}}{\text{d}v}g(v)\right|<\infty.
\label{eq:gain requirement}
\end{equation}
A wide class of electrical circuits can be described within these constructs. An example is Chua's circuit~\cite{Matsumoto-1985,Chua_Passivity99}, for which the impedance $z_\mathrm{osc}$ and nonlinear function $g(\cdot)$ are illustrated in Fig.~\ref{Fig:Examples}(a). The function $g(\cdot)$ is piecewise linear, and satisfies~\eqref{eq:gain requirement}. In previous work on voltage synchronization of voltage source inverters in small-scale power systems~\cite{Johnson_TCAS13,Johnson_JPV13,Johnson_TPELS13}, we introduced a nonlinear Li\'{e}nard-type dead-zone oscillator for which $z_\mathrm{osc}$ and $g(\cdot)$ are illustrated in Fig.~\ref{Fig:Examples}(b). In this case, the function $g(\cdot)$ is constructed with a negative resistance and a dead-zone function with finite slope, and thus satisfies~\eqref{eq:gain requirement}. Similar dead-zone-type oscillators have also been proposed in~\cite{Torres_Power12,Torres_Power13} for related power-systems applications. Some families of hyperchaotic circuits and negative-resistance oscillators can also be described with the model above, see\cite{Gismero-1990,Kyprianidis-1999,Malhame-2000,Kyprianidis-2001} and the references therein.  

A notable example of a well-known circuit that \emph{cannot} be described within the above framework is the Van der Pol oscillator~\cite{Khalil_Book02}. While the linear subsystem of the Van der Pol oscillator is the same as the nonlinear dead-zone oscillator, the nonlinear voltage-dependent current source, $g(v) \propto v^3$, which does not satisfy the slope requirement in~\eqref{eq:gain requirement} (see Fig.~\ref{Fig:Examples}(c)). 

\subsection{Electrical Network Model}\label{sub:networkmodel}
The nonlinear circuits are coupled through a passive, connected, LTI electrical network. The nodes of the electrical network are collected in the set $\mathcal{A}$, and branches of the electrical network are represented by the set of edges $\mathcal{E} := \{(m,n)\} \subset \mathcal{A} \times \mathcal{A}$. Let $\mathcal{N}:=\{1,\dots,N\} \subseteq \mathcal{A}$ collect {\em boundary nodes} that the nonlinear circuits are connected to, and let $\mathcal{I} = \mathcal{A} \setminus \mathcal{N}$ be the set of \emph{interior nodes} where the current injections are zero since they are not connected to the nonlinear electrical circuits. The series admittance corresponding to the $(m,n) \in \mathcal{E}$ branch is given by $y_{mn} \in \mathbb{C}$, and the shunt admittance connected between the $m^\mathrm{th}$ node and electrical ground is given by $y_{m} \in \mathbb{C}$. We will assume that the boundary nodes in the set $\mathcal{N}$ are not connected to any passive shunt elements, which implies $y_m = 0$ for all $m \in \mathcal{N}$. 

Denote the vectors that collect the nodal current injections and node voltages in the network by $i_\mathcal{A}$ and $v_\mathcal{A}$, respectively. The coupling between the circuits can be described by Kirchhoff's and Ohm's laws, which read in matrix-vector form as 
\begin{equation}
i_\mathcal{A}=Y_\mathcal{A}v_\mathcal{A}.\label{eq:i=Yv}
\end{equation}
In~\eqref{eq:i=Yv}, $Y_\mathcal{A} \in \mathbb{C}^{|\mathcal{A}|\times |\mathcal{A}|}$ denotes the admittance matrix of the electrical network, and it is constructed element-wise as 
\begin{equation}
[Y_\mathcal{A}]_{mn} := \left \{\begin{array}{ll} 
{y}_{m} + \sum_{(m,k) \in \mathcal{E}} y_{mk} , & \textrm{if } m=n,\\ 
- y_{mn} , &\textrm{if }(m,n)\in \mathcal{E},\\
0, & \textrm{otherwise}, \end{array} \right.
\label{eq:Ymatrix}
\end{equation}
where $y_{m}$ denotes the shunt admittance at node $m$ and $y_{mn} = y_{nm}$ denotes the line admittance of branch $(m,n)$. Notice that if the electrical network has no shunt elements, that is, $y_m = 0\, \,\forall m \in \mathcal{A}$, then $Y_\mathcal{A}$ is a singular matrix with zero row and column sums.

Let $i \!=\! [i_1,\dots,i_N]^\mathrm{T}$ and $v \!=\! [v_1,\dots,v_N]^\mathrm{T}$ be the vectors collecting the current injections and terminal voltages of the nonlinear circuits, and let $i_\mathcal{I}$ and $v_\mathcal{I}$ be the vectors collecting the current injections and nodal voltages for the interior nodes.\footnote{To be consistent with notation, we would have to include the subscript $\mathcal{N}$ when referring to the current and voltage vectors corresponding to the nonlinear circuits. However, we drop this subscript to ease exposition.} With this notation in place, we can rewrite~\eqref{eq:i=Yv} as 
\begin{equation}
\begin{bmatrix} i \\ i_\mathcal{I} \end{bmatrix} = \begin{bmatrix} Y_{\mathcal{N}\mathcal{N}} & Y_{\mathcal{N}\mathcal{I}} \\ Y_{\mathcal{N}\mathcal{I}}^\mathrm{T} & Y_{\mathcal{I}\mathcal{I}} \end{bmatrix} \begin{bmatrix} v \\ v_\mathcal{I} \end{bmatrix}.
\label{eq:OhmsLawMatrix}
\end{equation}  
Since the internal nodes are only connected to passive LTI circuit elements, all the entries of $i_\mathcal{I}$ are equal to zero in \eqref{eq:OhmsLawMatrix}. 

\begin{figure}[b]
\begin{centering}
\includegraphics{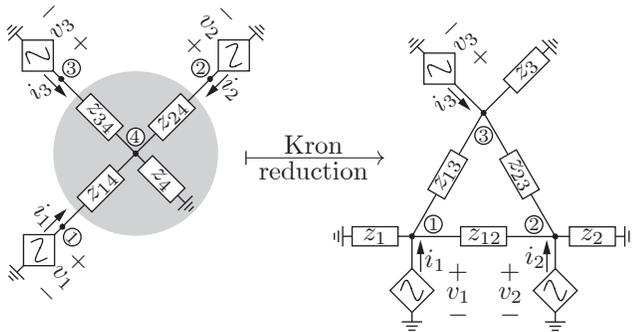}
\par\end{centering}
\caption{Kron reduction illustrated for a representative system comprising $N=3$ nonlinear circuits, where $\mathcal{A}=\{1,2,3,4\}$, $\mathcal{N}=\{1,2,3\}$, and $\mathcal{I}=\{4\}$. Original electrical network described by the admittance matrix, $Y_\mathcal{A}$, is shaded. In this setup Kron reduction is equivalent to the well-known star-delta transformation.\label{Fig:IllustratingKronReduction}}
\end{figure}

In the following, we assume that the submatrix $Y_\mathcal{II}$ is nonsingular such that the second set of equations in \eqref{eq:OhmsLawMatrix} can be uniquely solved for the interior voltages as  $v_\mathcal{I} = -Y_\mathcal{II}^{-1} Y_{\mathcal{N}\mathcal{I}}^\mathrm{T}v$. For $RL$, $RC$, or $RLC$ networks without shunt elements, the matrix $Y_{\mathcal A}$ is irreducibly block diagonally dominant (due to connectivity of the network), and $Y_{\mathcal{II}}$ is always nonsingular~\cite[Corollary~6.2.27]{RAH-CRJ:85}. For $RLC$ networks with shunt elements, it is possible to construct pathological cases where $Y_{\mathcal{II}}$ is singular, and the interior voltages $v_\mathcal{I}$ are not uniquely determined. In this paper, we assume that all principal submatrices are nonsingular and such pathological cases do not occur. Substituting $v_\mathcal{I} = -Y_\mathcal{II}^{-1} Y_{\mathcal{N}\mathcal{I}}^\mathrm{T}v$ in~\eqref{eq:OhmsLawMatrix}, then provides the following equations that relate the nonlinear-circuit current injections and terminal voltages:   
\begin{equation}
i = \left(Y_{\mathcal{NN}} - Y_{\mathcal{NI}} Y_{\mathcal{II}}^{-1} Y_{\mathcal{NI}}^\mathrm{T} \right) v =:Y v. \label{eq:OhmsLaw}
\end{equation}
This model reduction through a \emph{Schur complement} \cite{Zhang05} of the admittance matrix is known as \emph{Kron reduction}~\cite{Dorfler-13}. We refer to the matrix $Y$ in \eqref{eq:OhmsLaw} as the \emph{Kron-reduced admittance matrix}. From a control-theoretic perspective, \eqref{eq:OhmsLaw} is a {minimal realization} of the circuit \eqref{eq:OhmsLawMatrix}. We remark that even though the Kron-reduced admittance matrix $Y$ is well-defined, it is not necessarily the admittance matrix of a  passive circuit. Figure~\ref{Fig:IllustratingKronReduction} depicts an illustrative electrical network and its Kron-reduced counterpart for a system of $N=3$ nonlinear circuits.  

The results in this paper apply to Kron-reduced admittance matrices that satisfy the following two properties: 
\begin{enumerate}
\item[$\mathrm{(P1)}$] The Kron-reduced admittance matrix, $Y$, commutes with the projector matrix, $\Pi$, that is, $\Pi Y = Y \Pi$. 
\item[$\mathrm{(P2)}$] The Kron-reduced admittance matrix, $Y$, is normal, that is, $YY^* = Y^* Y$. Consequently, $Y$ can be diagonalized by a unitary matrix, that is, we can write $Y=Q \Lambda Q^*$, where $Q Q^* = I$ and $\Lambda$ is a diagonal matrix with diagonal entries composed of the eigenvalues of $Y$. 
\end{enumerate}
We will find $\mathrm{(P1)}$ useful in Section~\ref{sub:Compartmentalization}, where we derive a compartmentalized system description that clearly differentiates the linear and nonlinear subsystems in the electrical network. Similarly, $\mathrm{(P2)}$ will be leveraged in the proof of Theorem~\ref{Theorem: sync result I} in Section~\ref{sub:GAS}. We will identify classes of networks with and without shunt elements that satisfy properties $\mathrm{(P1)}$ and $\mathrm{(P2)}$ in Sections \ref{sec:GASSufficient1} and \ref{sec:GASSufficient2}, respectively. 

\section{Problem Statement and System Compartmentalization} \label{sec:Statement}
In this section, we first formulate the problem statement relating to global asymptotic synchronization. Next, we systematically compartmentalize the linear and nonlinear subsystems in the electrical network of coupled nonlinear electrical circuits. The differential $\mathcal{L}_2$ gains of the linear and nonlinear electrical subsystems will be used in subsequent sections to establish sufficient synchronization conditions. 

\subsection{Global Asymptotic Synchronization} \label{sub:GAS}
We are interested in global asymptotic synchronization of the terminal voltages of identical nonlinear circuits described in Section~\ref{sub:System-Description}, coupled through the electrical LTI network described in Section~\ref{sub:networkmodel}. In particular, we will seek sufficient conditions that ensure
\begin{equation}
\lim_{t\to\infty}v_{j}(t)-v_{k}(t)=0 \quad \forall j,k=1,\dots,N.\label{eq:GASCondition}
\end{equation}
For ease of analysis, we will find it useful to implement a coordinate transformation by employing the projector matrix, $\Pi$, to obtain the corresponding \emph{differential system} that clearly highlights signal differences. To emphasize the analytical advantages afforded by this coordinate transformation, note that: 
\begin{equation}
\widetilde{v}(t)^{\text{T}}\widetilde{v}(t)=\left(\Pi v(t)\right)^{\text{T}}\left(\Pi v(t)\right)=\frac{1}{2N}\sum_{j,k=1}^{N}\left(v_{j}(t)-v_{k}(t)\right)^{2}.
\end{equation}
Hence,~\eqref{eq:GASCondition} can be equivalently reformulated as
\begin{equation}
\lim_{t\to\infty}\widetilde{v}(t)=\lim_{t\to\infty}\Pi v(t) =0.
\end{equation}
The coordinate transformation with the projector matrix allows us to cast the voltage synchronization problem as an equivalent stability problem in the coordinates of the corresponding differential system. 

\subsection{Compartmentalization of Linear and Nonlinear Subsystems} \label{sub:Compartmentalization}
In order to establish synchronization conditions, we seek a system description where the linear and nonlinear subsystems ($z_\mathrm{osc}$ and $g(\cdot)$, respectively) in the network of coupled nonlinear circuits are clearly compartmentalized. In light of the importance of differential signals in facilitating the derivation of synchronization conditions, the compartmentalization is carried out in the corresponding differential system. 

\begin{figure}[b]
\begin{centering}
\includegraphics{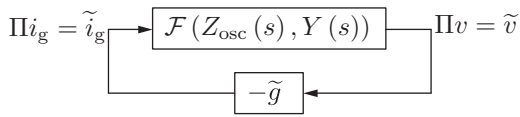}
\par\end{centering}
\caption{Block-diagram representation of the corresponding differential system. 
The linear and nonlinear portions of the system are compartmentalized in 
$\mathcal{F}\left(\cdot,\cdot\right)$ and $\widetilde{g}$, respectively.\label{Fig: ClosedLoopGeneralIncremental}}
\end{figure}

Towards this end, recall that the vectors $i$ and $v$ collect the current injections and terminal voltages of the nonlinear circuits. Similarly, denote by $i_\mathrm{g} := [i_{\mathrm{g}1}, \dots, i_{\mathrm{g}N}]^\mathrm{T}$, the vector that collects the currents sourced by the nonlinear voltage-dependent current sources. From Fig.~\ref{Fig:SingleOscillator}, we see that the terminal voltage of the $j^{\mathrm{th}}$ nonlinear circuit, $v_{j}$, can be expressed as
\begin{equation}
v_{j}=z_{\mathrm{osc}}\left(i_{\mathrm{g}j}-i_{j}\right),\,\forall j = 1,\dots,N.
\end{equation}
By collecting all $v_j$'s, we can write
\begin{equation}
v = Z_{\mathrm{osc}}\left(i_{\mathrm{g}}-i\right) =  Z_{\mathrm{osc}}i_{\mathrm{g}}-Z_{\mathrm{osc}}Yv,\label{eq:oscillator system v as function of currents}
\end{equation}
where $Z_{\mathrm{osc}}:=z_{\mathrm{osc}} I\in\mathbb{C}^{N\times N}$, and we substituted $i=Yv$ from~\eqref{eq:OhmsLaw}. A multiplication of both sides of~\eqref{eq:oscillator system v as function of currents} by the projector matrix $\Pi$ yields the differential terminal-voltage~vector
\begin{align}
\widetilde{v} = \Pi v & =\Pi\left(Z_{\mathrm{osc}}\left(i_{\mathrm{g}}-Yv\right)\right) =  Z_{\mathrm{osc}}\left(\Pi i_{\mathrm{g}}-\Pi Yv\right)\nonumber \\
 & =  Z_{\mathrm{osc}}\left(\widetilde{i}_{\mathrm{g}}-Y\widetilde{v}\right),\label{eq:diff sys step 1}
\end{align}
where we utilized the fact that $\Pi Z_{\mathrm{osc}}=\Pi z_{\mathrm{osc}}I=z_{\mathrm{osc}}I\Pi=Z_{\mathrm{osc}}\Pi$, and leveraged property $\mathrm{(P1)}$ which requires that the admittance and projector matrices commute. We can
now isolate $\widetilde{v}$ in~\eqref{eq:diff sys step 1} as follows:
\begin{equation}
\widetilde{v} = \left(I+Z_{\mathrm{osc}}Y\right)^{-1}Z_{\mathrm{osc}}\widetilde{i}_{\mathrm{g}}=  \mathcal{F}\left(Z_{\mathrm{osc}},Y\right)\widetilde{i}_{\mathrm{g}},\label{eq:differential isrc to v}
\end{equation}
where $\mathcal{F}\left(Z_{\mathrm{osc}},Y\right)$ is the \emph{linear fractional transformation} that captures the negative feedback interconnection of $Z_{\mathrm{osc}}$ and $Y$ (see~\eqref{eq:LFT} for a formal definition). Using~\eqref{eq:differential isrc to v}, we see that the corresponding differential system admits the compact block-diagram representation in Fig.~\ref{Fig: ClosedLoopGeneralIncremental}. The linear and nonlinear portions of the system are clearly compartmentalized by $\mathcal{F}\left(Z_{\mathrm{osc}},Y\right)$ and the map $\widetilde{g}:\widetilde{v}\to-\widetilde{i}_{\mathrm{g}}$, respectively. 
 
\section{Global Asymptotic Synchronization In Networks Without Shunt Elements} \label{sec:GASSufficient1}
This section focuses exclusively on synchronization in electrical networks that have no shunt elements. We begin this section by describing the class of electrical networks without shunt elements that satisfy $\mathrm{(P1)}$-$\mathrm{(P2)}$, and then present sufficient conditions for global asymptotic synchronization of the nonlinear-circuit terminal voltages in such networks.

First, we present a result which helps us to establish that Kron-reduced admittance matrices satisfy $\mathrm{(P1)}$ if the originating electrical networks have no shunt elements. The following result also offers a converse statement to \cite[Lemma 3.1]{Dorfler-13}.
\begin{thm} 
\label{Theorem: closure of zero row and column sums}
The following statements are equivalent:
\begin{enumerate}

	\item[(i)] The original electrical network has no shunt elements.

	\item[(ii)]The Kron-reduced network has no shunt elements.

\end{enumerate}
\end{thm}
\begin{IEEEproof} 
Let us first prove the sufficiency \emph{(i)}$\implies$\emph{(ii)}. In the absence of shunt elements, the admittance matrix $Y_\mathcal{A}$ has zero row sums by construction (see~\eqref{eq:Ymatrix}), that is,
\begin{equation}
\begin{bmatrix} \mathbf 0 \\ \mathbf 0 \end{bmatrix} = \begin{bmatrix} Y_{\mathcal{N}\mathcal{N}} & Y_{\mathcal{N}\mathcal{I}} \\ Y_{\mathcal{N}\mathcal{I}}^\mathrm{T} & Y_{\mathcal{I}\mathcal{I}} \end{bmatrix} \begin{bmatrix} \mathbf 1 \\ \mathbf 1 \end{bmatrix} \,.
\label{eq:ZeroRowSums}
\end{equation}  
An elimination of the second set of equations in \eqref{eq:ZeroRowSums} results in $\left(Y_{\mathcal{NN}} - Y_{\mathcal{NI}} Y_{\mathcal{II}}^{-1} Y_{\mathcal{NI}}^\mathrm{T} \right) \mathbf 1 = Y \mathbf 1 = \mathbf 0$, that is, $Y$ has zero row sums (and zero column sums due to closure of symmetry under the Schur complement~\cite{Zhang05}). By construction of the admittance matrix in~\eqref{eq:Ymatrix}, it follows that the Kron-reduced electrical network corresponding to $Y$ has no shunt elements.

\begin{figure*}[t]
\begin{center}
\subfigure[]{\includegraphics{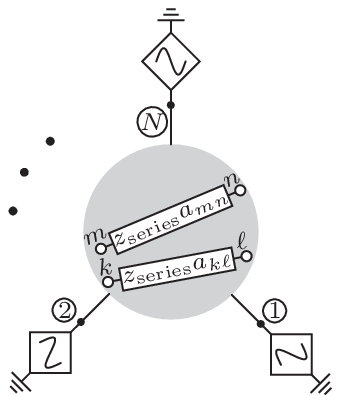}}\,
\subfigure[]{\includegraphics{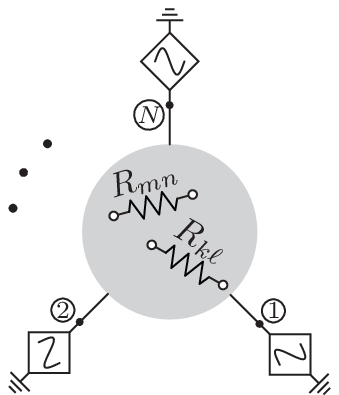}}\,
\subfigure[]{\includegraphics{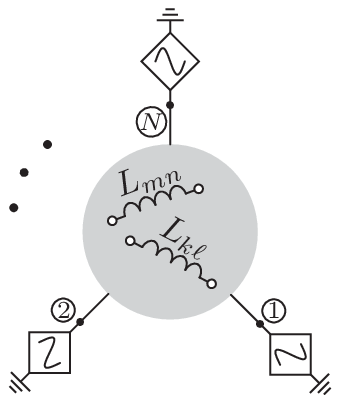}}\,
\subfigure[]{\includegraphics{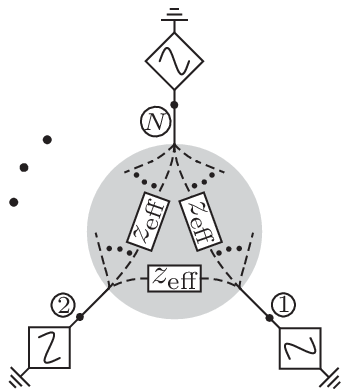}}
\end{center}
\caption{Illustrating electrical networks with no passive shunt elements that are: (a) uniform, (b) resistive, (c) lossless (inductive in this particular example), and (d) homogeneous. Uniform networks have identical per-unit-length impedances, $z_\mathrm{series}$, such that the value of the $(n,m)$ line impedance is expressed as the product of $z_\mathrm{series}$ and $a_{nm}$, the line length, which corresponds to the nonnegative weight from the underlying weighted Laplacian, $L$. In resistive (lossless) networks, notice that the line impedances are purely resistive (inductive in this particular example), and the values of the line impedances between different internal nodes are allowed to be arbitrary.  Finally, in homogeneous networks, the effective impedance between any two nodes is the same.}
\label{Fig:Networks}
\end{figure*}

We now prove the converse statement \emph{(ii)}$\implies$\emph{(i)} by proving its negation $\neg$\emph{(i)}$\implies$$\neg$\emph{(ii)}, that is, an original network with shunt elements always leads to a Kron-reduced network with shunt elements. Consider the augmented matrix $\widehat Y_{\mathcal A}$ associated with $Y_{\mathcal A}$, which is obtained by modeling the ground as an additional node in the network with index $|\mathcal A|+1$ and fixed (zero) voltage (see Appendix~\ref{AppA} for more details on the augmentation). Then $\widehat Y_{\mathcal A}$ is the admittance matrix associated to a network without shunt elements, that is, $\widehat Y_{\mathcal A} \mathbf 1=\mathbf 0$. By the reasoning (i)$\implies$(ii) above, the associated Kron-reduced admittance matrix $\widehat Y$ has no shunt elements, that is, $\widehat Y \mathbf 1=\mathbf 0$. Since Kron-reduction and augmentation commute (see Lemma~\ref{Lemma: Commutativity of Augmentation and Kron reduction} in Appendix~\ref{AppA}), the reduced network with admittance matrix $Y$ obtained by removing the grounded node (that is removing the column and row with index $N+1$ from $\widehat Y$) does not have zero row sums. Equivalently, the Kron-reduced network has  shunt elements.\end{IEEEproof}

\begin{cor}
If the original electrical network has no shunt elements, then property $\mathrm{(P1)}$ holds, that is, $Y$ commutes\,with~$\Pi$.
\end{cor}
\begin{IEEEproof} 
If the original network has no shunt elements, then the Kron-reduced network has no shunt elements. Thus, $Y$ has zero row and column sums, and $Y$ commutes with $\Pi$.
\end{IEEEproof}

\subsection{Identifying Electrical Networks that Satisfy $\mathrm{(P1)}$-$\mathrm{(P2)}$}\label{sub:NetworkTypes1}
The synchronization criteria we develop within this section apply to the following classes of networks: 
\begin{enumerate}
	\item[(i)] {\em networks with uniform line characteristics} \cite{SYC-PT:12}, in which the branch admittances are given $y_{nm} = \subscr{y}{series}a_{nm}$ for all $ (n,m) \in \mathcal{E}$, where $a_{nm} \in \mathbb R$ is real-valued and $\subscr{y}{series} \in \mathbb C \setminus \{0\}$ is identical for every branch (see Figs.~\ref{Fig:Networks}(a)-(c)); 
	\item[(ii)] {\em homogeneous networks} \cite{Dorfler-13}, in which the effective impedances are identical for all boundary nodes, that is, $z_{nm} =: z_\mathrm{eff} =  r + \mathrm{j} x, r, x \in \mathbb{R}, \forall n,m \in \mathcal{N}$ (see Fig.~\ref{Fig:Networks}(d)).
\end{enumerate}
For these networks, we will first derive the general form of the Kron-reduced admittance matrices, and then subsequently demonstrate compliance to $\mathrm{(P2)}$. 

We first focus on networks with uniform line characteristics, which physically correspond to networks for which all branches are made of the same material \cite{SYC-PT:12}, i.e., the admittance of each branch $(n,m)$ depends on its constant per-unit-length admittance, $\subscr{y}{series} \in \mathbb C$, and its length, $a_{nm}>0$ (see Fig.~\ref{Fig:Networks}(a)). Notice that these networks include as special cases resistive networks (Fig.~\ref{Fig:Networks}(b)) and lossless networks (Fig.~\ref{Fig:Networks}(c)) for which $\subscr{y}{series}$ is real-valued or purely imaginary, respectively. For these networks, we can express $Y_{\mathcal A} = \subscr{y}{series} \cdot L_{\mathcal A}$, where $L_{\mathcal A}$ is a symmetric, positive semidefinite, and real-valued Laplacian matrix. We have the following result:

\begin{lem}
\label{Lemma: closure of uniformity}
Consider a network with uniform line characteristics, that is, $Y_{\mathcal A} = \subscr{y}{series} \cdot L_{\mathcal A}$, where $L_{\mathcal A} \in \mathbb R^{|\mathcal A|\times|\mathcal A|}$ is a real-valued Laplacian matrix and $\subscr{y}{series} \in \mathbb C$. Then, the Kron-reduced network has uniform line characteristics with the Kron-reduced admittance matrix given by 
\begin{equation}
	Y =  \subscr{y}{series} L 
	\label{eq:YUniform}
\end{equation}
where $L= L_{\mathcal{NN}} - L_{\mathcal{NI}} L_{\mathcal{II}}^{-1} L_{\mathcal{NI}}^\mathrm{T}$.
\end{lem}
The proof of Lemma \ref{Lemma: closure of uniformity} follows by direct construction of the Kron-reduced matrix and due to the closure properties of Kron reduction~\cite[Lemma 3.1]{Dorfler-13}. Due to the special form of the Kron-reduced admittance matrix in~\eqref{eq:YUniform}, it follows that $Y$ is diagonalizable with a unitary matrix. We conclude that $\mathrm{(P2)}$ is satisfied by Kron-reduced matrices for which the original electrical network has uniform line characteristics.

To address homogeneous networks, we recall from \cite[Theorem III.4]{Dorfler-13} that (in the purely resistive case) a sparse electrical network becomes denser under Kron reduction and even complete under mild connectivity assumptions. However, the branch admittances in the reduced network are still heterogeneous and reflect the topology and electrical properties of the original network. In the following result, we show that for a homogeneous original network, the associated Kron-reduced network is characterized by identical branch admittances.

\begin{lem}
\label{Lemma: homogeneous loopless graphs}
The following statements are equivalent:
\begin{enumerate}

	\item[(i)] The original network is homogeneous:  for all boundary nodes $n,m \in \{1,\dots,N\}$, the pairwise effective impedances take the uniform value $z_{nm} = \subscr{z}{eff} \in \mathbb C \setminus \{0\}$.
		
	\item[(ii)] The Kron-reduced network is complete and the branch admittances take the uniform value $\subscr{y}{series} \in \mathbb C \setminus \{0\}$. Equivalently, the Kron-reduced admittance matrix is given by
\begin{equation}
	Y  = \subscr{y}{series} \Gamma\,,
	\label{eq:YHomogeneous}
\end{equation}
where $\Gamma = N I - \mathbf{1}\mathbf{1}^\mathrm{T}$ is the Laplacian of the complete graph.
\end{enumerate}
If statements (i) and (ii) are true, then $\subscr{z}{eff} = \frac{2}{N \subscr{y}{series}}$.
\end{lem}

Lemma~\ref{Lemma: homogeneous loopless graphs} is obtained as a direct corollary to Theorem~\ref{Theorem: Kron Reduction in Homogeneous Networks} in Appendix~\ref{AppA}. Since the Laplacian of the complete graph is $\Gamma=N\Pi$, it commutes with the projector matrix $\Pi$. Finally, notice that since the Kron-reduced admittance matrix $Y$ in the homogeneous case~\eqref{eq:YHomogeneous} is a special case of the Kron-reduced admittance matrix for the network with uniform line characteristics~\eqref{eq:YUniform}, it follows that Kron-reduced admittance matrices for homogeneous electrical networks satisfy $\mathrm{(P2)}$. 

\subsection{Sufficient Condition for Global Asymptotic Synchronization} \label{sub:GAS}
This subsection derives sufficient conditions to ensure global asymptotic synchronization in the network of coupled nonlinear circuits described in Section~\ref{sub:NetworkTypes1}. First, we present a lemma that establishes an upper bound on the differential $\mathcal{L}_2$ gain of the function $g(\cdot)$, that governs the nonlinear voltage-dependent current sources in the nonlinear circuits.

\begin{lem}(\cite[Lemma 1]{Johnson_TCAS13})
The differential $\mathcal{L}_{2}$ gain of $g(\cdot)$ is finite, and upper bounded by $\sigma$:
\begin{equation}
\widetilde{\gamma}\left(g\right):=\frac{\Vert \widetilde{i}_{\mathrm{g}}\Vert _{\mathcal{L}_{2}}}{\left\Vert \widetilde{v}\right\Vert _{\mathcal{L}_{2}}}\leq \sigma := \sup_{v\in\mathbb{R}}\left|\frac{\mathrm{d}}{\mathrm{d}v}g(v)\right|<\infty.\label{eq:Lemma 1 result}
\end{equation}
\end{lem}
We now provide a sufficient synchronization condition for the case where the nonlinear circuits are connected in networks with uniform line characteristics. Subsequently, we consider homogeneous networks. 
\begin{thm}
\label{Theorem: sync result I}
Suppose the electrical network that couples the system of $N$ identical nonlinear circuits has no shunt elements, and has uniform line characteristics. Let the Kron-reduced admittance matrix be of the form $Y = y_\mathrm{series}  L$ as in \eqref{eq:YUniform}. The terminal voltages of the nonlinear circuits synchronize in the sense of~\eqref{eq:GASCondition} if for all $j \in\{ 2,\dots,N\}$
\begin{equation}
\Vert \mathcal{F}(z_\mathrm{osc}(\mathrm{j}\omega),y_\mathrm{series}(\mathrm{j}\omega) \lambda_j)\Vert_\infty \sigma < 1 \,,
\label{eq:GAS}
\end{equation}
where $\lambda_j$, $j \in\{ 2,\dots,N\}$, are the nonzero eigenvalues of the Laplacian matrix $L$. 
\end{thm}
For purely resistive networks the synchronization condition~\eqref{eq:GAS} has to be evaluated only for $\lambda_2$. The second-smallest eigenvalue $\lambda_{2}$ of the Laplacian matrix is known as the {\em algebraic connectivity}, and it is a spectral connectivity measure~\cite{MF:73}. It can be shown that the algebraic connectivity in a resistive Kron-reduced network upper-bounds the algebraic connectivity in the original network  \cite[Theorem III.5]{Dorfler-13}. Hence, condition~\eqref{eq:GAS} implies that the nonlinear circuits should be sufficiently strongly connected, which is aligned with synchronization results in complex oscillator networks with a static interconnection topology (i.e., with only resistive elements)~\cite{Dorfler-13-Synch}. On the other hand, if the interconnecting network is {\em dynamic}, e.g., if it contains capacitive or inductive storage elements, then the  synchronization condition \eqref{eq:GAS} needs to be evaluated for all nonzero network modes $\lambda_j$, $j \in \{2,\dots, N\}$.

\begin{IEEEproof}[Proof of Theorem \ref{Theorem: sync result I}]
Consider the block-diagram of the differential system in Fig.~\ref{Fig: ClosedLoopGeneralIncremental}. From Lemma 1, we have
\begin{equation}
\Vert \widetilde{i}_\mathrm{g}\Vert _{\mathcal{L}_{2}}\leq\sigma\Vert \widetilde{v}\Vert _{\mathcal{L}_{2}}.
\label{eq:lemma1}
\end{equation}
For the linear fractional transformation, we can write 
\begin{equation}
\Vert \widetilde{v}\Vert _{\mathcal{L}_{2}}\leq\widetilde{\gamma}\left(\mathcal{F}\left(Z_{\mathrm{osc}},Y\right)\right)\Vert \widetilde{i}_\mathrm{g}\Vert _{\mathcal{L}_{2}}+\eta,\label{eq:proof step 1}
\end{equation}
for some non-negative $\eta$, where $\widetilde{\gamma}\left(\mathcal{F}\left(Z_{\mathrm{osc}},Y\right)\right)$ denotes the differential $\mathcal{L}_{2}$ gain of the linear fractional transformation. By combining~\eqref{eq:lemma1} and~\eqref{eq:proof step 1}, we arrive at
\begin{equation}
\Vert \widetilde{v}\Vert _{\mathcal{L}_{2}}\leq\widetilde{\gamma}\left(\mathcal{F}\left(Z_{\mathrm{osc}},Y\right)\right)\sigma\Vert \widetilde{v}\Vert _{\mathcal{L}_{2}}+\eta.\label{eq:L2_vdiff1}
\end{equation}
By isolating $\Vert \widetilde{v}\Vert _{\mathcal{L}_{2}}$ from~\eqref{eq:L2_vdiff1}, we can write 
\begin{equation}
\Vert \widetilde{v}\Vert _{\mathcal{L}_{2}}\leq\frac{\eta}{1-\widetilde{\gamma}\left(\mathcal{F}\left(Z_{\mathrm{osc}},Y\right)\right)\sigma},
\end{equation}
provided that the following condition holds
\begin{equation}
\widetilde{\gamma}\left(\mathcal{F}\left(Z_{\mathrm{osc}},Y\right)\right)\sigma<1.\label{eq:General Synch Result}
\end{equation}
If \eqref{eq:General Synch Result} holds true, then we have $\widetilde{v}\in\mathcal{L}_{2}$. It follows from Barbalat's lemma \cite{Hamadeh_Constructive10,Hamadeh_Designing12,Hamadeh_Thesis} 
that $\lim_{t\to\infty}\widetilde{v}(t)=\mathbf{0}$. Hence, if the network of nonlinear circuits satisfies the condition~\eqref{eq:General Synch Result}, global asymptotic synchronization can be guaranteed. 

In the remainder of the proof, we establish an equivalent condition for~\eqref{eq:General Synch Result}. By definition of the differential $\mathcal{L}_{2}$ gain of the linear fractional transformation, we can express
\begin{align}
& \widetilde{\gamma}\left(\mathcal{F}\left(Z_{\text{osc}},Y\right)\right)  = \widetilde{\gamma}\left(\mathcal{F}\left(z_\mathrm{osc}I,Y\right)\right) \\
&= \underset{\omega\in\mathbb{R}}{\mathrm{sup}}\,\frac{\left\Vert \mathcal{F}\left(z_\mathrm{osc}\left(\mathrm{j}\omega\right)I,Y\left(\mathrm{j}\omega\right)\right)\widetilde{i}_{\mathrm{g}}\left(\mathrm{j}\omega\right)\right\Vert _{2}}{\left\Vert \widetilde{i}_{\mathrm{g}}(\mathrm{j}\omega)\right\Vert _{2}}\nonumber \\
& =\underset{\omega\in\mathbb{R}}{\mathrm{sup}}\,\frac{\left\Vert \left(I+z_\mathrm{osc}(\mathrm{j}\omega)Y\left(\mathrm{j}\omega\right)\right)^{-1}z_\mathrm{osc}(\mathrm{j}\omega)\widetilde{i}_{\mathrm{g}}\left(\mathrm{j}\omega\right)\right\Vert _{2}}{\left\Vert \widetilde{i}_{\mathrm{g}}(\mathrm{j}\omega)\right\Vert _{2}}\nonumber \\
& =\underset{\omega\in\mathbb{R}}{\mathrm{sup}}\,\frac{\left\Vert Q\left(I+z_\mathrm{osc} (\mathrm{j}\omega) \subscr{y}{series}(\mathrm{j}\omega)  \Lambda\right)^{-1}z_\mathrm{osc}(\mathrm{j}\omega)Q^\mathrm{T}\widetilde{i}_{\mathrm{g}}(\mathrm{j}\omega)\right\Vert _{2}}{\left\Vert Q^\mathrm{T}\widetilde{i}_{\mathrm{g}}(\mathrm{j}\omega)\right\Vert _{2}}, \nonumber \label{eq:proof step 2}
\end{align}
where we made use of property $\mathrm{(P2)}$ to diagonalize the admittance matrix as $Y=\subscr{y}{series} L = \subscr{y}{series} Q\Lambda Q^\mathrm{T}$, where $Q$ is unitary and $\Lambda$ is a diagonal matrix containing the real-valued and nonnegative eigenvalues of the Laplacian matrix\,$L$.

Since the Kron-reduced network has no shunt elements connected to ground, the row and column sums of $Y$ are zero. Furthermore, since the Kron-reduced network is connected, $Y$ has a single zero eigenvalue. Analogous comments apply to $L$ and we obtain $\lambda_1=0$ with corresponding eigenvector $q_{1}=(1/\sqrt{N}) \mathbf{1}$. Finally, since $\mathbf{1}^{\mathrm{T}}\Pi=\mathbf{0}^{\mathrm{T}}$, we can express
\begin{equation}
Q^\mathrm{T}\widetilde{i}_{\mathrm{g}}=Q^\mathrm{T}\Pi i_{\mathrm{g}}=\left[0,\, p\right]^{\mathrm{T}},\label{eq:observation 1}
\end{equation}
where  $p\in\mathbb{C}^{N-1}$ is made of the non-zero elements of $Q^\mathrm{T}\Pi i_{\mathrm{g}}$. 

Using the observation in~\eqref{eq:observation 1}, denoting the $N-1\times N-1$ identity matrix by $I_{N-1}$, and the diagonal matrix with entries composed of the non-zero eigenvalues of $Y$ by $\Lambda_{N-1}$, we can simplify~\eqref{eq:proof step 1} as follows:
\begin{align}
&\widetilde{\gamma}\left(\mathcal{F}\left(z_\mathrm{osc}I,Y\right)\right)\nonumber \\
 & =\underset{\omega\in\mathbb{R}}{\mathrm{sup}}\,\frac{\left\Vert \left(I_{N-1}+z_\mathrm{osc}(\mathrm{j}\omega)\subscr{y}{series}(\mathrm{j}\omega) \Lambda_{N-1}\right)^{-1}z_\mathrm{osc}(\mathrm{j}\omega)p(\mathrm{j}\omega)\right\Vert _{2}}{\left\Vert p(\mathrm{j}\omega)\right\Vert _{2}}\nonumber \\
 & =\max_{j=2,\dots,N}\,\underset{\omega\in\mathbb{R}}{\mathrm{sup}}\,\left|\frac{z_\mathrm{osc}(\mathrm{j}\omega) }{1+z_\mathrm{osc}(\mathrm{j}\omega)\subscr{y}{series}(\mathrm{j}\omega) \lambda_j}\right| \,.
\label{eq:proof part 4}
\end{align}
By combining~\eqref{eq:proof part 4} and~\eqref{eq:General Synch Result}, and using the definition of the linear fractional transformation in~\eqref{eq:LFT}, we arrive at condition~\eqref{eq:GAS}.~\end{IEEEproof}

\begin{cor} 
\label{Corollary: sync result II}
Suppose the original electrical network that couples the system of $N$ identical nonlinear circuits has no shunt elements and is homogeneous with a Kron-reduced admittance matrix $Y = y_\mathrm{series} \Gamma$ as in~\eqref{eq:YHomogeneous}. The terminal voltages of the nonlinear circuits connected in such a network synchronize in the sense of~\eqref{eq:GASCondition} if
\begin{equation}
\Vert \mathcal{F}(z_\mathrm{osc}(\mathrm{j}\omega),y_\mathrm{series}(\mathrm{j}\omega) N)\Vert_\infty \sigma < 1. 
\label{eq:GAS1}
\end{equation}
\end{cor}

\begin{IEEEproof}
The condition in~\eqref{eq:GAS1} follows from the fact that the eigenvalues of $\Gamma$ (the Laplacian matrix of the complete graph), are given by $\lambda_1 = 0$ and $\lambda_2=\dots=\lambda_N=N$.
\end{IEEEproof}

\section{Global Asymptotic Synchronization In Networks With Shunt Elements} \label{sec:GASSufficient2}
In this section, we explore the family of electrical networks with shunt elements for which sufficient synchronization conditions similar to~\eqref{eq:GAS} can be derived. 

\subsection{Identifying Electrical Networks that Satisfy $\mathrm{(P1)}$-$\mathrm{(P2)}$}\label{sub:NetworkTypes2}
Consider the case of a Kron-reduced electrical network, where---in addition to a single nonlinear circuit---each node $ m \in \mathcal{N}$ in the network is connected to an identical shunt admittance $y_m = y_\mathrm{shunt}$. In this case, the Kron-reduced admittance matrix can be expressed as 
\begin{equation}
\label{eq:AdmittanceSpecial1}
Y = y_\mathrm{shunt}I + Y',
\end{equation}
where $Y'$ corresponds to the admittance matrix that captures the coupling between the nonlinear circuits. Note that the row and column sums of $Y'$ are zero since $y'_{m} = 0$ for all $m \in \mathcal{N}$ by construction (since $Y'$ does not include shunt elements). If the network modeled by such a $Y'$ has uniform line characteristics (such as resistive or lossless lines), then we obtain
\begin{equation}
	\label{eq:AdmittanceSpecial2}
	Y = y_\mathrm{shunt}I + y_\mathrm{series} L, 
\end{equation}
where $L$ is the associated real-valued, symmetric Laplacian matrix. Clearly, $Y$ as in \eqref{eq:AdmittanceSpecial2} satisfies properties $\mathrm{(P1)}$ and $\mathrm{(P2)}$.

In general, it is difficult to identify networks that admit Kron-reduced admittance matrices of the form in~\eqref{eq:AdmittanceSpecial1} or even \eqref{eq:AdmittanceSpecial2}. However, we can identify a family of electrical networks that admit Kron-reduced admittance matrices of the same form as~\eqref{eq:AdmittanceSpecial2}. Towards this end, we first present a result on Kron reduction of homogeneous networks with shunt elements. 

\begin{lem}
\label{Lemma: homogeneous loopy graphs}
The following statements are equivalent:
\begin{enumerate}

	\item[(i)] The original network is homogeneous: for all boundary nodes $n,m \in \{1,\dots,N\}$, the pairwise effective impedances take the uniform value $z_{nm} = \subscr{z}{eff--series} \in \mathbb C \setminus \{0\}$ and the effective impedances to electrical ground (denoted by the node, $N+1$) take the uniform value $z_{n(N+1)} = \subscr{z}{eff--shunt} \in \mathbb C \setminus \{0\}$, with $\frac{\subscr{z}{eff--series}}{\subscr{z}{eff--shunt}} \neq \frac{2 N}{N-1}$.
		
	\item[(ii)] The branch and shunt admittances in the Kron-reduced network are uniform, that is, there are $\subscr{y}{series} \in \mathbb C \setminus \{0\}$ and $\subscr{y}{shunt} \in \mathbb C \setminus \{0\}$, with $\subscr{y}{shunt} \neq - N \subscr{y}{series}$, such that
	\begin{equation}
Y = y_\mathrm{shunt}I + y_\mathrm{series} \Gamma.
\label{eq:YredParticular}
\end{equation}

\end{enumerate}
If statements (i) and (ii) are true, then 
\begin{align}
	\subscr{y}{series} &= \frac{2}{2N \subscr{z}{eff--shunt} - (N-1) \subscr{z}{eff--series}}
	\,, \nonumber
	\\
	\subscr{y}{shunt} &= \frac{2(\subscr{z}{eff--series} - 2\subscr{z}{eff--shunt})}{\subscr{z}{eff--series}((N-1)\subscr{z}{eff--series} - 2 N \subscr{z}{eff--shunt})}
	\,.
\label{eq:Correspondence}
\end{align}  
\end{lem}

Lemma~\ref{Lemma: homogeneous loopy graphs} is obtained as a direct corollary to Theorem~\ref{Theorem: Kron Reduction in Homogeneous Networks} in Appendix~\ref{AppA}, by inverting the constitutive relations between admittances and effective impedances. The parametric assumption $\subscr{y}{shunt} \neq - N \subscr{y}{series}$ (respectively, ${\subscr{z}{eff--series}}/{\subscr{z}{eff--shunt}} \neq {2 N}/({N-1})$) is practically not restrictive: it is violated only in pathological cases, e.g., when a capacitive (respectively inductive) shunt load compensates exactly for $N$ inductive (respectively capacitive) line flows.

\begin{figure}[t]
\begin{centering}
\includegraphics{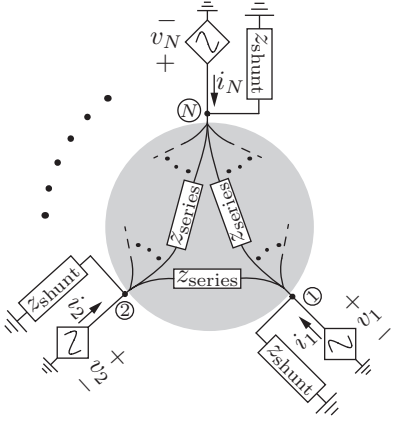}
\par\end{centering}
\caption{Kron-reduced electrical network recovered from a homogeneous originating network. The shaded region captures the inter-circuit interactions through identical line impedances that are equal to $z_\mathrm{series}$. All the shunt impedances are equal to $z_\mathrm{shunt}$.}
\label{fig:HomogeneousKron}
\end{figure}

The admittance matrix $Y$ in~\eqref{eq:YredParticular} satisfies $\mathrm{(P1)}$ and $\mathrm{(P2)}$, and it is clearly a special case of~\eqref{eq:AdmittanceSpecial1}, with $Y' = y_\mathrm{series} \Gamma$. While the formulation in~\eqref{eq:YredParticular} is more restrictive, Lemma~\ref{Lemma: homogeneous loopy graphs} identifies the {\em unique} class of electrical networks that admit Kron-reduced matrices of the form~\eqref{eq:YredParticular}. An illustration of a Kron-reduced electrical network recovered from a homogeneous originating electrical network is depicted in Fig.~\ref{fig:HomogeneousKron}.

\subsection{Sufficient Condition for Global Asymptotic Synchronization} \label{sub:GAS1}

We now present sufficient conditions for global asymptotic synchronization for the cases where the Kron-reduced admittance matrices are given by~\eqref{eq:AdmittanceSpecial2}, or as a special case, by~\eqref{eq:YredParticular}. 

\begin{cor} 
\label{Corollary: sync result III}
Suppose the electrical network that couples the system of $N$ identical nonlinear circuits admits a Kron-reduced admittance matrix given by~\eqref{eq:AdmittanceSpecial2}, where the network corresponding to $L$ has uniform line characteristics and no shunt elements connected to ground. The terminal voltages of the nonlinear circuits synchronize in the sense of~\eqref{eq:GASCondition} if for all $j \in\{ 2,\dots,N\}$
\begin{equation}
\Vert \mathcal{F}(z_\mathrm{eq}(\mathrm{j}\omega),\subscr{y}{series}(\mathrm{j}\omega)\lambda_j)\Vert_\infty \sigma < 1 \,,
\label{eq:GAS2}
\end{equation}
where $\lambda_j$, $j \in\{ 2,\dots,N\}$, are the nonzero eigenvalues of the Laplacian matrix $L$, and $z_\mathrm{eq}$ is the equivalent impedance of the parallel combination of $z_\mathrm{shunt}:=y_\mathrm{shunt}^{-1}$ and $z_\mathrm{osc}$ given by
\begin{equation}
z_\mathrm{eq} := \frac{z_\mathrm{shunt} z_\mathrm{osc}}{z_\mathrm{shunt} + z_\mathrm{osc}}.
\label{eq:zeta}
\end{equation}
\end{cor}

\begin{IEEEproof}
The proof for this corollary follows along the same lines as that for Theorem~\ref{Theorem: sync result I}. In particular, if~\eqref{eq:General Synch Result} holds, then synchronization is guaranteed. Now, consider that with $Y = y_\mathrm{shunt}I + \subscr{y}{series} L$, we have that
\begin{align}
&\mathcal{F}\left(Z_{\mathrm{osc}},Y\right)=\left(I+Z_{\mathrm{osc}}Y\right)^{-1}Z_{\mathrm{osc}}\nonumber \\
&=\left(I+Z_{\mathrm{osc}}\left(y_\mathrm{shunt}I+\subscr{y}{series} L\right)\right)^{-1}Z_{\mathrm{osc}}\nonumber \\
&=\left(I+\frac{z_\mathrm{shunt} z_{\mathrm{osc}}}{z_\mathrm{shunt}+z_{\mathrm{osc}}}\subscr{y}{series} L\right)^{-1}\frac{z_\mathrm{shunt} z_{\mathrm{osc}}}{z_\mathrm{shunt}+z_{\mathrm{osc}}}I \nonumber \\
&=\mathcal{F}\left(z_\mathrm{eq} I,\subscr{y}{series} L\right),\label{eq:cor step 1}
\end{align}
where the last line in~\eqref{eq:cor step 1} follows from the definition of the linear fractional transformation~\eqref{eq:LFT}, and the definition of $z_\mathrm{eq}$ in~\eqref{eq:zeta}. 
By repeating the reasoning as in the proof of Theorem~\ref{Theorem: sync result I}, we obtain
\begin{align}
&\widetilde{\gamma}\left(\mathcal{F}\left(Z_\mathrm{osc},Y\right)\right)=\widetilde{\gamma}\left(\mathcal{F}\left(z_\mathrm{eq} I,\subscr{y}{series}L\right)\right)\nonumber \\
& =\max_{j=2,\dots,N}\,\,\underset{\omega\in\mathbb{R}}{\mathrm{sup}}\, \,\, \left|\frac{z_\mathrm{eq}(\mathrm{j}\omega)}{1+z_\mathrm{eq}(\mathrm{j}\omega)\subscr{y}{series}(\mathrm{j}\omega)\lambda_j}\right|.
\label{eq:proof part 6}
\end{align}
The claimed synchronization condition \eqref{eq:GAS2} then follows by combining~\eqref{eq:General Synch Result} and \eqref{eq:proof part 6}. 
\end{IEEEproof}

\begin{cor} 
\label{Corollary: sync result IV}
Suppose the original electrical network that couples the system of $N$ identical nonlinear circuits admits a Kron-reduced admittance matrix of the form~\eqref{eq:YredParticular}. The terminal voltages of nonlinear circuits connected in such a network synchronize in the sense of~\eqref{eq:GASCondition} if
\begin{equation}
\Vert \mathcal{F}(z_\mathrm{eq}(\mathrm{j}\omega),y_\mathrm{series}(\mathrm{j}\omega) N)\Vert_\infty \sigma < 1. 
\label{eq:GAS4}
\end{equation}
\end{cor}

\begin{IEEEproof}
The condition in~\eqref{eq:GAS4} follows from the fact that the eigenvalues of $\Gamma$ (the Laplacian matrix of the complete graph), are given by $\lambda_1 = 0$ and $\lambda_2=\dots=\lambda_N=N$.
\end{IEEEproof}

\begin{figure}[t]
\begin{centering}
\includegraphics{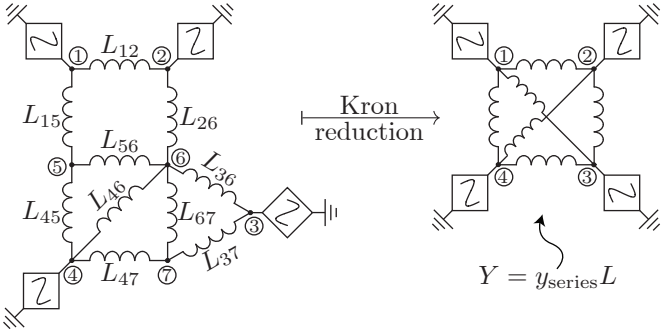}
\par\end{centering}
\caption{Schematic of lossless inductive network and the Kron-reduced counterpart examined in Section~\ref{sec:CS1}. This is an example of a network with uniform line characteristics (see Section~\ref{sec:GASSufficient1}), in that the Kron-reduced admittance matrix can be expressed as $Y = y_\mathrm{series} L$, where $L$ is a weighted, real-valued, symmetric Laplacian matrix.}
\label{fig:Lossless}
\end{figure}

\section{Case Studies} \label{sec:CaseStudies}
We now present simulation case studies to validate the synchronization conditions in some illustrative LTI electrical network topologies that interconnect {\em Chua's circuits}~\cite{Matsumoto-1985}. The electrical schematic of Chua's circuit is depicted in Figure~\ref{Fig:Examples}(a). Two network topologies are considered in this section: i) an arbitrary lossless network with nonidentical inductive line impedances and without shunt elements, and ii) a homogeneous network with a symmetric star topology and a shunt load. 

The nonlinear voltage-dependent current source in Chua's circuit is illustrated in Fig.~\ref{Fig:Examples}(a), and the impedance of the linear subsystem is given by
\begin{equation}
z_{\mathrm{osc}}(s)=\frac{RLC_\mathrm{a}C_\mathrm{b} s^3 + LC_\mathrm{a} s^2 + RC_\mathrm{a} s}{RLC_\mathrm{a}C_\mathrm{b} s^3 + (C_\mathrm{a}^2+LC_\mathrm{b}+LC_\mathrm{a})s^2 + RC_\mathrm{a} s + 1 }.
\label{eq:ChuaImpedance}
\end{equation}
Since the function $g(\cdot)$ in Fig.~\ref{Fig:Examples}(a) is piecewise linear, it immediately satisfies~\eqref{eq:gain requirement}. Parameters of the constituent circuit elements in Chua's circuit utilized in both case studies are summarized in Appendix~\ref{sub:ModelParameters}. With the choice of circuit parameters, it follows that $\sigma = \sup_{v\in\mathbb{R}}\left|\frac{\text{d}}{\text{d}v}g(v)\right| = \sigma_2$.

\begin{figure}[t!]
\vspace{-0.2in}
\centering
{ \subfigure[]{\includegraphics{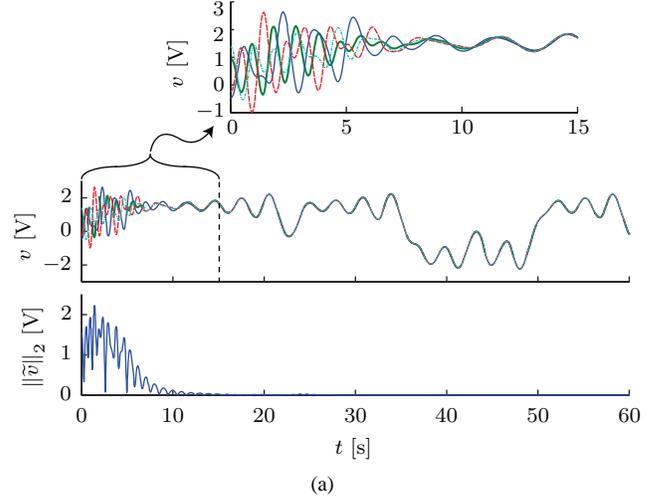}}\vspace{-0.05in}
   \subfigure[]{\includegraphics{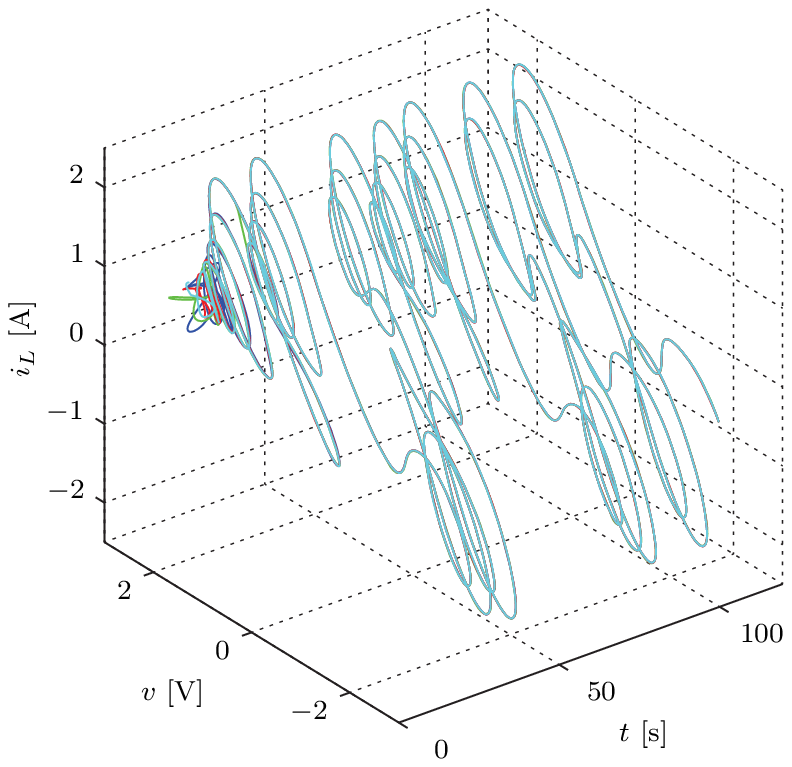}}\vspace{-0.05in}
\caption{Synchronization of terminal voltages in Chua's circuits for the lossless electrical network depicted in Fig.~\ref{fig:Lossless}. (a) Terminal voltages, $v(t)$, and voltage synchronization error, $||\widetilde{v}(t)||_2$. The inset illustrates a close-up view of system dynamics at startup with nonidentical initial conditions. (b) Chaotic double-scroll attractor is discernible in the asymptotic limit.}
\label{fig:LnetSynch}
}
\end{figure}

\begin{figure}[t!]
\begin{centering}
\includegraphics{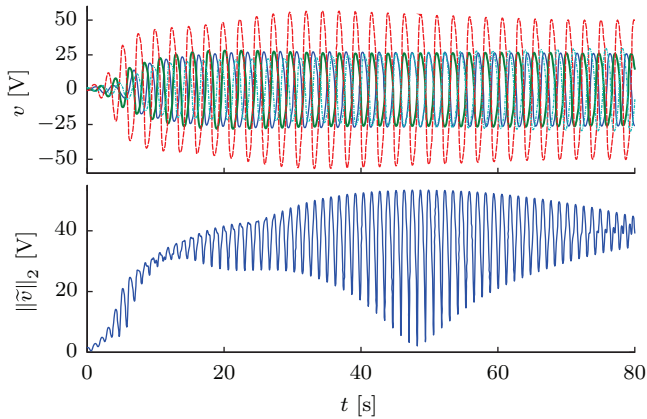}
\par\end{centering}
\caption{Terminal voltages, $v(t)$, and voltage synchronization error,
  $||\widetilde{v}(t)||_2$, for a particular network where
  $||\mathcal{F}\left(z_{\mathrm{osc}}(\mathrm{j}\omega),y_{\mathrm{series}}^{-1}(\mathrm{j}\omega\right)
  \lambda_j || _{\infty} \sigma \nless 1$, $j=2, 3, 4$, and synchronization is
  not guaranteed.}
\label{fig:LnetNoSynch}
\end{figure}

\subsection{Lossless Inductive Network without Shunt Elements} \label{sec:CS1}
The network topology examined here is illustrated in Fig.~\ref{fig:Lossless}. This is an example of a network with uniform line characteristics (see Section~\ref{sub:NetworkTypes1}). Following Lemma~\ref{Lemma: closure of uniformity}, we obtain the Kron-reduced network (also illustrated in Fig.~\ref{fig:Lossless}) with admittance matrix given by $Y = y_\mathrm{series} L$, where $L$ is a weighted, real-valued, symmetric Laplacian matrix. The sufficient synchronization condition for  this case is given by~\eqref{eq:GAS}.

For the first set of network parameters in Appendix~\ref{sub:NetworkParametersLossless}, we get $||\mathcal{F}\left(z_{\mathrm{osc}}(\mathrm{j}\omega),y_{\mathrm{series}}^{-1}(\mathrm{j}\omega\right) \lambda_j || _{\infty} \sigma<1$, $j = 2, 3, 4$, which implies that the terminal voltages of the Chua's circuits are guaranteed to synchronize. We confirm this with time-domain simulations. Figure~\ref{fig:LnetSynch}(a) illustrates the terminal voltages and the voltage synchronization error, with the inset capturing a close-up during startup with nonidentical initial conditions as the voltages begin to pull into phase. Figure~\ref{fig:LnetSynch}(b) depicts a three-dimensional view of the internal states of the Chua's circuits as a function of time, and clearly demonstrates the chaotic double-scroll attractor \cite{Matsumoto-1985} in the asymptotic limit.

Now consider the second set of network parameters in Appendix~\ref{sub:NetworkParametersLossless}. For this set of parameters, we obtain $||\mathcal{F}\left(z_{\mathrm{osc}}(\mathrm{j}\omega),y_{\mathrm{series}}^{-1}(\mathrm{j}\omega\right) \lambda_j|| _{\infty} \sigma \nless 1$, $j = 2, 3, 4$, and our condition \eqref{eq:GAS} is violated. While this is not an indication that the terminal voltages cannot synchronize (since the condition \eqref{eq:GAS} is only sufficient), it turns out that in this case, the terminal voltages indeed do not synchronize. With the same set of initial conditions as before, we plot the the terminal voltages and the voltage synchronization error in Fig.~\ref{fig:LnetNoSynch} for this network. 

\begin{figure}[b]
\begin{centering}
\includegraphics{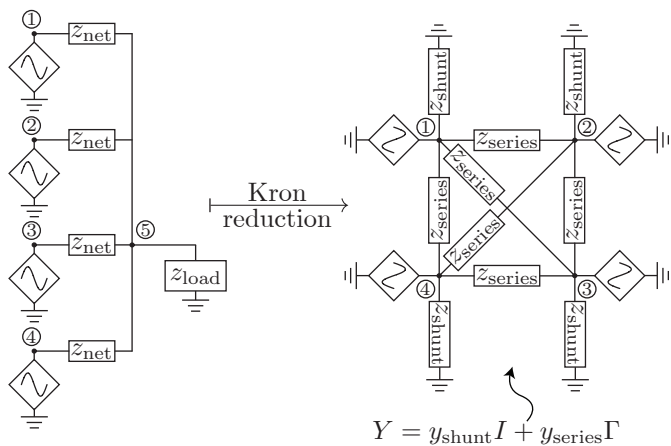}
\par\end{centering}
\caption{Network composed of nonlinear electrical circuit connected to a common load through identical branch impedances, and the corresponding Kron-reduced circuit. Homogeneity of the original electrical network implies that the admittance matrix of the Kron-reduced network is given by $Y = y_\mathrm{shunt}I + y_\mathrm{series} \Gamma$.}
\label{fig:CaseStudyBNetwork}
\end{figure}

\subsection{Homogeneous Network with Shunt Load}
\begin{figure}[t]
\begin{centering}
\includegraphics{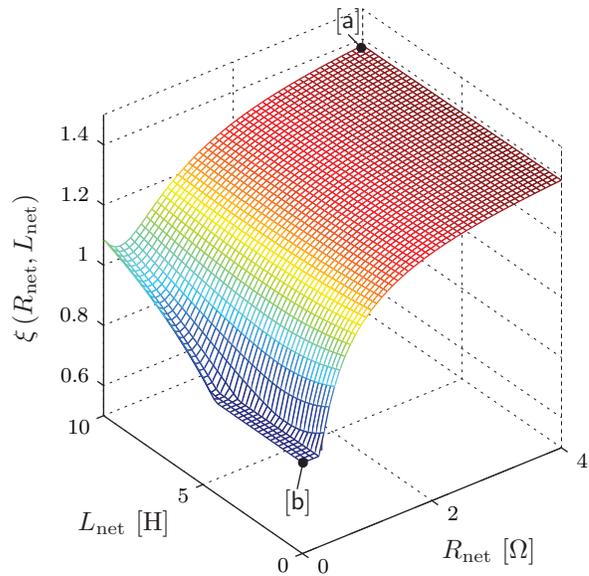}
\par\end{centering}
\caption{The function, $\xi(R_\mathrm{net},L_\mathrm{net})$ plotted for a range of values of $R_\mathrm{net}$ and $L_\mathrm{net}$. Synchronization is guaranteed for values of $R_\mathrm{net}$ and $L_\mathrm{net}$ where $\xi(R_\mathrm{net},L_\mathrm{net})<1$.}
\label{fig:threedee}
\end{figure}

The network topology examined here is illustrated in Fig.~\ref{fig:CaseStudyBNetwork}. The network branch impedance is given by $z_{\mathrm{net}}\left(s\right)=sL_{\mathrm{net}}+R_{\mathrm{net}}$. The impedance of the shunt load connected to the internal node is denoted by $z_{\mathrm{load}}\left(s\right)$. This is an example of a homogeneous network since the effective impedances between any two nonlinear circuits are identical, and the effective impedances between the nonlinear circuits and electrical ground are also identical (see Section~\ref{sub:NetworkTypes2}). In particular
\begin{equation}
z_\mathrm{eff-series} = 2z_\mathrm{net}, \quad z_\mathrm{eff-shunt} = z_\mathrm{net} + z_\mathrm{load}.
\end{equation}
Following Lemma~\ref{Lemma: homogeneous loopy graphs}, we obtain the Kron-reduced network (also illustrated in Fig.~\ref{fig:CaseStudyBNetwork}) with admittance matrix given by $Y = y_\mathrm{shunt}I + y_\mathrm{series} \Gamma$, where $\Gamma$ is the Laplacian of the complete graph. By applying~\eqref{eq:Correspondence}, we obtain
\begin{align}
y_\mathrm{shunt} &= \left(z_{\mathrm{net}} + 4 z_{\mathrm{load}}\right)^{-1}, \nonumber \\
y_\mathrm{series} &= z_{\mathrm{load}} \left(z_{\mathrm{net}} \left(z_{\mathrm{net}} + 4 z_{\mathrm{load}}\right)^{-1}\right)^{-1}.
\end{align}
Substituting $y_\mathrm{shunt}$ and $y_\mathrm{series}$ in~\eqref{eq:GAS4}, we obtain the following sufficient condition for synchronization in this network:
\begin{equation}
\left\Vert \mathcal{F}(z_\mathrm{osc},z_\mathrm{net}^{-1}) \right \Vert_\infty \sigma = \left\Vert \frac{z_\mathrm{osc}(\mathrm{j}\omega) z_\mathrm{net}(\mathrm{j}\omega)}{z_\mathrm{osc}(\mathrm{j}\omega) + z_\mathrm{net}(\mathrm{j}\omega)} \right \Vert_\infty \sigma < 1.
\label{eq:SynchConditionParallelCaseStudy}
\end{equation}
It is instructive to explore the impact of the network parameters $R_\mathrm{net}$ and $L_\mathrm{net}$ on synchronization. First, consider Fig.~\ref{fig:threedee}, which plots $\xi(R_\mathrm{net},L_\mathrm{net}):=\Vert \mathcal{F}(z_\mathrm{osc},z_\mathrm{net}^{-1}) \Vert_\infty \sigma$ for a range of different values of $R_\mathrm{net}$ and $L_\mathrm{net}$. For sufficiently low values of $R_\mathrm{net}$ and $L_\mathrm{net}$, synchronization is guaranteed since $\xi(R_\mathrm{net},L_\mathrm{net})<1$ as required by~\eqref{eq:SynchConditionParallelCaseStudy}. 

\begin{figure}[t]
\begin{centering}
\includegraphics{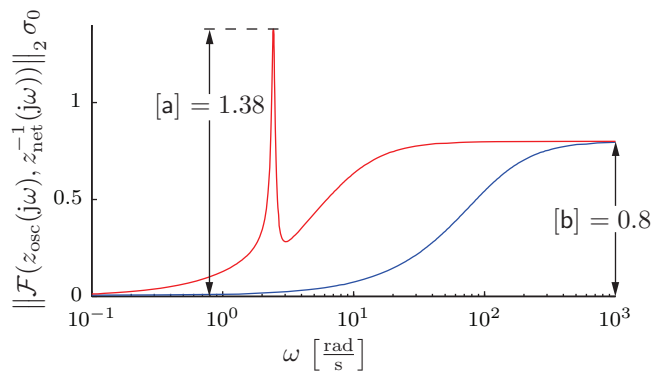}
\par\end{centering}
\caption{Magnitude of $\mathcal{F}(z_\mathrm{osc}(\mathrm{j} \omega),z_\mathrm{net}^{-1}(\mathrm{j} \omega)) \sigma$ as a function of frequency, $\omega$.}
\label{fig:XiBode}
\end{figure}

Consider now, the asymptotic behavior of $\xi(R_\mathrm{net},L_\mathrm{net})$. Particularly, we will focus on two points \textsf{[a]} and \textsf{[b]} that are marked in Fig.~\ref{fig:threedee}. Figure~\ref{fig:XiBode} depicts the magnitude of $\mathcal{F}(z_\mathrm{osc}(\mathrm{j} \omega),z_\mathrm{net}^{-1}(\mathrm{j} \omega)) \sigma$ as a function of frequency, $\omega$, for three representative values of $R_\mathrm{net}, L_\mathrm{net}$; two of which correspond to the asymptotes \textsf{[a]} and \textsf{[b]}. The effect of reducing the values of $R_\mathrm{net}, L_\mathrm{net}$ translates to damping the peak of the magnitude response. Synchronization is guaranteed when the peak is less than unity.   

\section{Concluding Remarks and Directions for Future Work} \label{sec:Conclusions}

We derived a synchronization condition for a class of nonlinear electrical circuits coupled through passive LTI electrical networks. We considered particular classes of networks, where perfect synchronization of the terminal voltages can be achieved. These classes included homogeneous networks and networks with uniform line characteristics---both with and without shunt elements. Whereas these classes of networks seem to be restrictive at first, it is the belief of the authors that---with the present setup---perfect synchronization cannot be achieved for more general and heterogeneous networks, where the  nonlinear circuits are possibly non-identical and support different loads. In this case, the subsystems have no common asymptotic dynamics to synchronize on~\cite{Wieland:2013tj}.

In ongoing and future work, we plan to address the problems of synchronization in heterogeneous networks and the regulation of the asymptotic synchronized dynamics. Further topics of interest include the analysis of Kron reduction of general $RLC$ circuits (including pathological cases) and synchronization through directed and possibly nonlinear electrical networks with diodes and rectifiers.

\begin{appendix}

\subsection{Kron Reduction of Complex-Symmetric Matrices} \label{AppA}
In this appendix, we discuss the Kron reduction of complex-valued admittance matrices and the properties of the effective impedances. The following results are extensions from the real-valued and symmetric cases considered in~ \cite{Dorfler-13}, to complex-symmetric (and not necessarily Hermitian) settings relevant in this work. Since only a subset of results in \cite{Dorfler-13} directly carries over to the complex-symmetric case, we present self-contained statements with brief proof sketches. These results directly lead up to Lemmas \ref{Lemma: homogeneous loopless graphs} and \ref{Lemma: homogeneous loopy graphs} in this paper.

First, notice that an admittance matrix $Y_{\mathcal A}$ without shunt elements is singular due to zero row and column sums, and an admittance matrix with at least one shunt element is invertible due to irreducibly block diagonally dominance \cite[Corollary~6.2.27]{RAH-CRJ:85}. To analyze regular and singular admittance matrices simultaneously, we associate an {\em augmented admittance matrix} $\widehat Y_{\mathcal A}$ to a regular admittance matrix $Y_{\mathcal A}$:
\begin{equation}
	\widehat Y_{\mathcal A}
	:=
	\left[\begin{array}{ccc|c}
	& & & - y_{1} \\
	& \mbox{\huge$Y_{\mathcal A}$} & & \vdots \\
	& & & - y_{|\mathcal A|} \\\hline
	- y_{1} & \cdots & -y_{|\mathcal A|} & \Bigr.\sum_{m=1}^{|\mathcal A|} y_{m}
	\end{array}\right]
	\,,
\end{equation}
The augmented admittance matrix $\widehat Y_{\mathcal A}$ corresponds to the case when the ground is modeled as an additional node with index\,$|\mathcal A|+1$ and zero voltage. Notice that $\widehat Y_{\mathcal A}$ is singular with zero row and column sums. Likewise, a singular admittance matrix resulting from a network without shunt elements can be regularized by grounding an arbitrary node. We denote the Kron-reduced matrix associated to $\widehat Y_{\mathcal A}$ by $\widehat Y$. As it turns out, the augmentation process and the Kron reduction commute. 

\begin{lem}
\label{Lemma: Commutativity of Augmentation and Kron reduction}
Consider a singular admittance matrix $Y_{\mathcal A}$, its augmented matrix $\widehat Y_{\mathcal A}$, and their associated Kron-reduced matrices $Y$ and $\widehat Y$, respectively. The following diagram commutes:
\begin{center}
	\includegraphics[width=0.6\columnwidth]{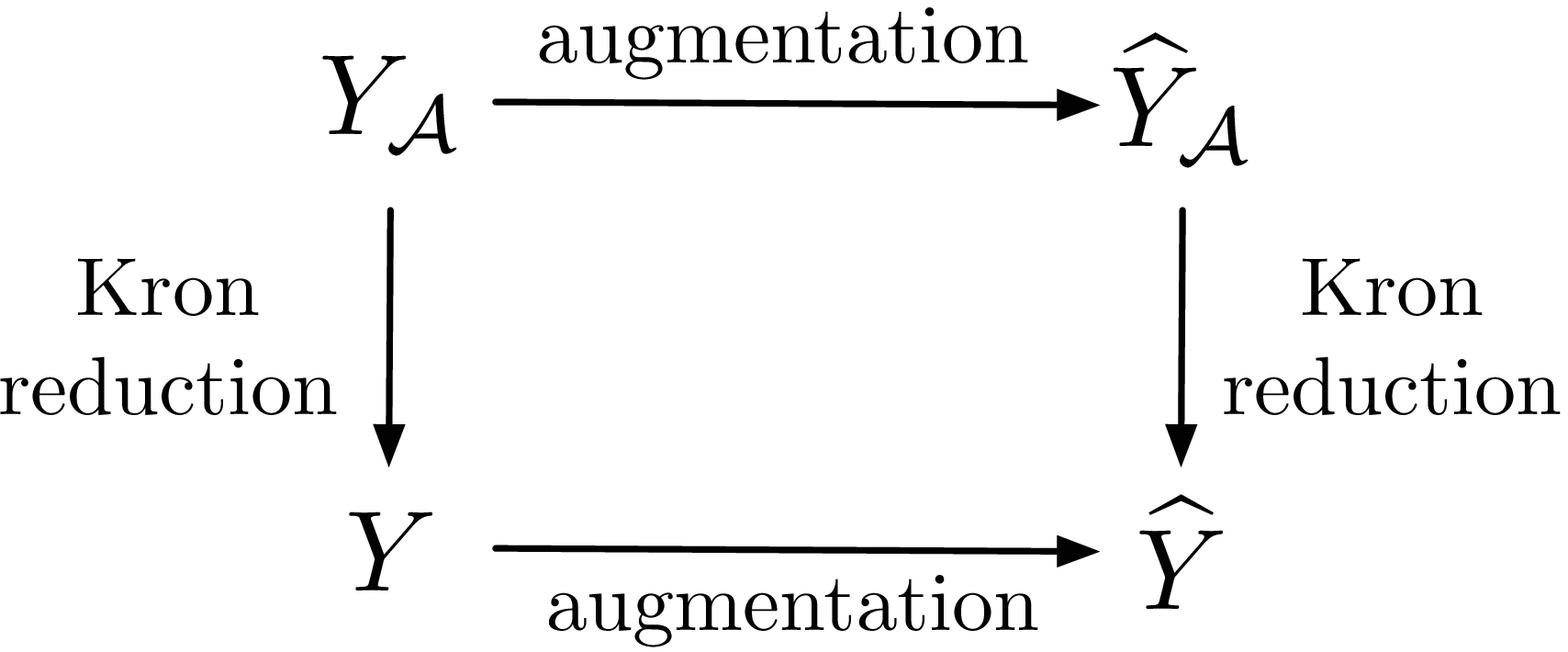}
\end{center}
\end{lem}

\begin{IEEEproof} 
The proof is analogous to the proof of \cite[Lemma III.1, Property 3]{Dorfler-13}. The result in \cite{Dorfler-13} relies on the Quotient Formula \cite[Theorem 1.4]{Zhang05} which extends to complex-valued matrices, as well as the closure of symmetry and zero row (column) sums under Kron reduction (shown in Theorem \ref{Theorem: closure of zero row and column sums}).~\end{IEEEproof}

As the next key property, we establish that the effective impedances $z_{nm}$ among the boundary nodes $n,m \in \mathcal N$ are invariant under Kron reduction and augmentation.

\begin{thm}
\label{Theorem: Invariance of Effective Impedance under Kron Reduction}
Consider the admittance matrix $Y_{\mathcal A}$ and the Kron-reduced matrix $Y$. The following statements hold:

\begin{enumerate}
	\item[1)] Invariance under Kron reduction: the effective impedance between any two boundary nodes is equal when computed from $Y$ or $Y_{\mathcal A}$, that is, for any $n,m \in \{1,\dots,N\}$ 
\begin{align}
	\!\!\!\!\!\!\!\!\!
	z_{nm} &=  (e_{n} - e_{m})^\mathrm{T} Y^{\dagger} (e_{n} - e_{m})  \nonumber \\
&\equiv (e_{n} - e_{m})^\mathrm{T} Y_{\mathcal A}^{\dagger} (e_{n} - e_{m}).
\end{align}
	\item[2)] Invariance under augmentation: if $Y_{\mathcal A}$ is a nonsingular matrix, then the effective impedance is equal when computed from $Y_{\mathcal A}$ or $\widehat Y_{\mathcal A}$, that is, for any $n,m \in \{1,\dots,|\mathcal A|\}$
\begin{align}
	\!\!\!\!\!\!\!\!\!\!
	z_{nm}	&= 	(e_{n}-e_{m})^\mathrm{T} Y_{\mathcal A}^{-1} (e_{n}-e_{m}) \nonumber \\
	& \equiv (e_{n}-e_{m})^\mathrm{T} \widehat Y_{\mathcal A}^{\dagger} 	(e_{n}-e_{m}).
\end{align}
\end{enumerate}
Equivalently, statements 1) and 2) imply that, if $Y_{\mathcal A}$ is a regular admittance matrix, then the following diagram commutes:
\begin{center}
	\includegraphics[width=0.7\columnwidth]{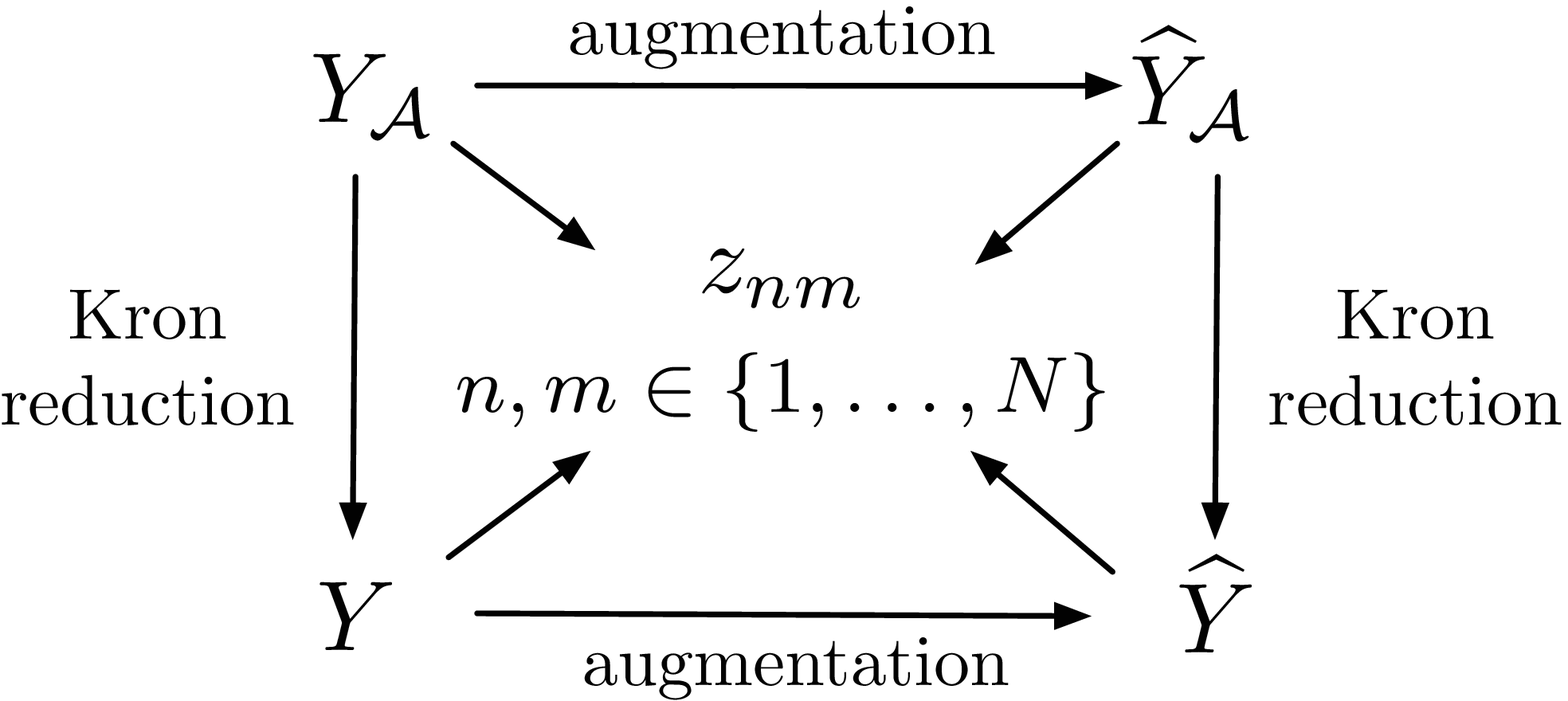}
\end{center}
\end{thm}

\begin{IEEEproof} 
To prove Theorem \ref{Theorem: Invariance of Effective Impedance under Kron Reduction}, we first establish some matrix identities. We need the following identity for a singular admittance matrix $Y \in \mathbb C^{N \times N}$ and a real nonzero number $\delta$:
\begin{equation}
	\left(Y + (\delta/N) \, \mathbf 1 \mathbf 1^{T} \right)^{-1}
	=
	Y^{\dagger} + (1/\delta N) \, \mathbf 1 \mathbf 1^{T} 
	\label{eq: matrix identity with pseudo inverse}
	\,.
\end{equation}
Using the projector formula (for a singular admittance matrix) $Y Y^{\dagger} = Y^{\dagger} Y = \Pi$, the identity \eqref{eq: matrix identity with pseudo inverse} can be verified since the product of the left-hand and the right-hand side of \eqref{eq: matrix identity with pseudo inverse} equal the identity matrix. If a singular admittance matrix $Y \in \mathbb C^{N \times N}$ is of dimension $N \geq 3$, then by taking the $N^\mathrm{th}$ node as a reference and deleting the associated $N$th column and $N^\mathrm{th}$ row, we obtain the nonsingular matrix $\overline Y  \in \mathbb C^{(N-1) \times (N-1)}$. As suggested by physical intuition, the effective impedance among the nodes $n,m \in \{1,\dots,N-1\}$ is not affected by grounding the $N^\mathrm{th}$ node, that is, for all $n,m \in \{1,\dots,N-1\}$
\begin{align}
	z_{nm} &= (e_{n} - e_{m})^\mathrm{T} \overline Y^{-1} (e_{n} - e_{m}) \nonumber \\
& \equiv (e_{n} - e_{m})^\mathrm{T} Y^{\dagger} (e_{n} - e_{m}).
	\label{eq: effective resistance from Dirichlet Laplacian}
\end{align}
The identity \eqref{eq: effective resistance from Dirichlet Laplacian} can be verified by using the formula
$
\overline Y^{-1}_{nm} 
=
Y^{\dagger}_{nm} - Y^{\dagger}_{nN} - Y^{\dagger}_{mN} + Y^{\dagger}_{NN}
$
\cite[Appendix B, formula (17)]{FF-AP-JMR-MS:07} whose derivation extends to the complex-valued case.

To prove statement 1), consider first the case when $Y_{\mathcal A}$ is invertible due to presence of shunt admittances. Recall that we are interested in the effective impedances only among the boundary nodes, that is, the leading principal $(N \times N)$-block of $Y_{\mathcal A}^{\dagger} = Y_{\mathcal A}^{-1}$. The Schur complement formula \cite[Theorem 1.2]{Zhang05} gives the leading $(N \times N)$-block of $Y_{\mathcal A}^{-1}$ as the inverse Schur complement $Y^{-1}$. It follows that, for all $n,m \in \{1,\dots,N\}$,
\begin{align}
	z_{nm} &=  (e_{n} - e_{m})^\mathrm{T} Y^{-1} (e_{n} - e_{m}) \nonumber \\ 
& \equiv (e_{n} - e_{m})^\mathrm{T} Y_{\mathcal A}^{-1} (e_{n} - e_{m}).
\end{align}
On the other hand, if $Y_{\mathcal A}$ is singular due to absence of shunt admittances, an analogous reasoning applies on the image of $Y$ and using the identity \eqref{eq: matrix identity with pseudo inverse}, or after grounding an arbitrary interior node (i.e., regularizing $Y_{\mathcal A}$) and using the identity \eqref{eq: effective resistance from Dirichlet Laplacian}, see the proof of \cite[Theorem III.8, Property 1]{Dorfler-13} for details.

To prove statement 2), notice that the regular admittance matrix $Y_{\mathcal A}$ with shunt elements is the leading principal $(|\mathcal A|\times|\mathcal A|)$-block of the augmented singular admittance matrix $\widehat Y_{\mathcal A}$. Statement 2) follows then directly from identity \eqref{eq: effective resistance from Dirichlet Laplacian} after replacing $N$, $Y$, and $\overline Y$ with $N+1$, $\widehat Y_{\mathcal A}$, and $Y_{\mathcal A}$.
\end{IEEEproof}

\begin{thm}
\label{Theorem: Kron Reduction in Homogeneous Networks}
Consider an admittance matrix $Y_{\mathcal A}$ and its Kron-reduced matrix $Y$. Consider the following two cases:

\smallskip
{1) No shunt elements:} Assume that $Y_{\mathcal A}$ is singular due to the absence of shunt elements. Let $Z \in \mathbb C^{N \times N}$ be the matrix of effective impedances. The following statements are equivalent:
\begin{enumerate}

	\item[(i)] The effective impedances among the boundary nodes $\{1,\dots,N\}$ are uniform, that is, there is $\subscr{z}{eff} \in \mathbb C \setminus \{0\}$ such that $z_{nm} = \subscr{z}{eff}$ for all distinct $n,m \in \{1,\dots,N\}$; 		
	\item[(ii)] The branch admittances in the Kron-reduced network take the uniform value $\subscr{y}{series} \in \mathbb C \setminus \{0\}$, that is, $Y = \subscr{y}{series} \Gamma$.

\end{enumerate}
If statements (i) and (ii) are true, then $\subscr{z}{eff} = \frac{2}{N \subscr{y}{series}}$.
\smallskip

{2) Shunt elements:} Assume that $Y_{\mathcal A}$ is regular due to the presence of shunt elements. Consider the grounded node $|\mathcal A|+1$ and the augmented admittance matrices $\widehat Y_{\mathcal A}$ and $\widehat Y$. Let $Z \in \mathbb R^{(|\mathcal A|+1) \times (|\mathcal A|+1)}$ be the matrix of effective impedances in the augmented network. The following statements are equivalent:
\begin{enumerate}

	\item[(iii)] The effective impedances both among the boundary nodes $\{1,\dots,N\}$ and between all boundary nodes $\{1,\dots, N\}$ and the grounded node $|\mathcal A|+1$ are uniform, that is, there are $\subscr{z}{eff--series} \in \mathbb C \setminus \{0\}$ and $\subscr{z}{eff--shunt} \in \mathbb C \setminus \{0\}$ satisfying ${\subscr{z}{eff--series}}/{\subscr{z}{eff--shunt}} \neq {2 N}/{N-1}$ such that $z_{ij} = \subscr{z}{eff--series}$ for all distinct $n,m \in \{1,\dots,N\}$ and $z_{n,|\mathcal A|+1} = \subscr{z}{eff--shunt}$ for all $n \in \{1,\dots,N\}$.
	
	\item[(iv)] The branch and shunt admittances in the Kron-reduced network are uniform, that is, there are $\subscr{y}{series} \in \mathbb C \setminus \{0\}$ and $\subscr{y}{shunt} \in \mathbb C \setminus \{0\}$ satisfying $\subscr{y}{shunt} \neq - N \subscr{y}{series}$ such that $Y = y_\mathrm{shunt}I + y_\mathrm{series} \Gamma$.

\end{enumerate}
If statements (iii) and (iv) are true, then 
\begin{align}
	\subscr{z}{eff--series} &= \frac{2}{N \subscr{y}{series} + \subscr{y}{shunt}} \,,
	\\
	\subscr{z}{eff--shunt} &= \frac{\subscr{y}{shunt}+\subscr{y}{series}}{\subscr{y}{shunt}(N \subscr{y}{series}+\subscr{y}{shunt})} \,.
\end{align}
\end{thm}

The admittance assumption $\subscr{y}{shunt} \neq - N \subscr{y}{series}$ (and the equivalent assumption ${\subscr{z}{eff--series}}/{\subscr{z}{eff--shunt}} \neq {2 N}/{N-1}$ for the effective impedances) guarantees regularity (respectively a single zero eigenvalue) of the admittance matrix. As discussed in Section \ref{sub:NetworkTypes2}, this assumption is practically not restrictive.

\begin{IEEEproof} [Proof of Theorem \ref{Theorem: Kron Reduction in Homogeneous Networks}]
Here, we present the proof strategy for case 1). Due to the invariance of the effective impedance (and Kron reduction) under the augmentation process (see Lemma \ref{Lemma: Commutativity of Augmentation and Kron reduction} and Theorem \ref{Theorem: Invariance of Effective Impedance under Kron Reduction}), an analogous reasoning and similar formulae apply to case 2), see \cite{Dorfler-13}. 

We first prove the statement $(i) \implies (ii)$: Notice that $Z=Z^\mathrm{T}$ has zero diagonal elements and $Y^{\dagger}$ is symmetric with zero row and column sums. Hence, both matrices have $N (N-1)/2$ independent elements, and the linear formula relating the elements $z_{nm}$ and $Y^{\dagger}_{nm}$ can be inverted \cite[identity (34)]{Dorfler-13}: 
	\begin{equation}
		Y^{\dagger}_{nm} 
		=
		-\frac{1}{2} \Bigl( z_{nm} - \frac{1}{N} \sum_{k=1}^{N} (z_{nk} + z_{mk}) 
		+ \frac{1}{N^{2}} \sum_{k, \ell=1}^{N} z_{k\ell} \Bigr)\,.
	\end{equation}
From the above formula, it can be readily verified that a uniform effective impedance matrix, $Z = \subscr{z}{eff} \, (\mathbf 1 \mathbf 1^\mathrm{T} - I)$, 
yields a uniform inverse matrix $Y^{\dagger} = \subscr{z}{eff}/(2N) \, \Gamma$. It is worth mentioning that for an admittance matrix with uniform branch admittances, $Y = \subscr{y}{series}\, \Gamma$, the pseudo inverse $Y^{\dagger}$ is again an admittance matrix with uniform branch admittances given by 
 \begin{equation}
 	Y^{\dagger} = \left( \subscr{y}{series}\, \Gamma \right)^{\dagger} =  1/(N^{2} \subscr{y}{series}) \, \Gamma \,.
	\label{eq: inverse of uniform Laplacian}
 \end{equation}
Identity \eqref{eq: inverse of uniform Laplacian} can be verified since $Y$ and $Y^{\dagger}$ satisfy the Penrose equations.  According to~\eqref{eq: inverse of uniform Laplacian}, this uniform inverse matrix $Y^{\dagger}=\subscr{z}{eff}/(2N) \, \Gamma$ yields the uniform admittance matrix $Y = {2}/({N \subscr{z}{eff}})\,\Gamma$. 

Now we consider the converse implication $(ii) \implies (i)$. Due to Theorem \ref{Theorem: Invariance of Effective Impedance under Kron Reduction} the effective impedance is invariant under Kron reduction. In this case, substituting $Y^{\dagger}$ from~\eqref{eq: inverse of uniform Laplacian} in~\eqref{eq:zeff}, we see that the effective impedances are given by
\begin{align}
z_{nm} &= \frac{1}{(N^{2} \subscr{y}{series}) }(e_n - e_m)^\mathrm{T} (N I - \mathbf{1} \mathbf{1}^\mathrm{T}) (e_n - e_m) \nonumber \\
& = \frac{2}{N \subscr{y}{series}}= \subscr{z}{eff}, \forall n,m \in \{1,\dots,N\}.
\end{align}
Hence, the $N (N-1)/2$ pairwise effective impedances $z_{nm}$ between the boundary nodes are uniform.\end{IEEEproof}

\subsection{Parameters of Chua's Circuits} \label{sub:ModelParameters}
Linear subsystem parameters: $R=10/7 \, \mathrm{\Omega}$, $L=1/7 \, \mathrm{H}$, $C_\mathrm{a}=1/9 \, \mathrm{F}$, $C_\mathrm{b}=1 \, \mathrm{F}$. Nonlinear-subsystem parameters: $\sigma_0=-0.8 \, \mathrm{S}$, $\sigma_1=-0.5 \, \mathrm{S}$, $\sigma_2=0.8 \, \mathrm{S}$, $\varphi_0=1 \, \mathrm{V}$, $\varphi_1=14 \, \mathrm{V}$.

\subsection{Lossless network parameters} \label{sub:NetworkParametersLossless}
\subsubsection{Guaranteed Synchronization}
$L_{12} = 0.834 \, \mathrm{H}$, $L_{15} = 0.671 \, \mathrm{H}$,
$L_{26} = 0.277 \, \mathrm{H}$, $L_{56} = 1.0575 \, \mathrm{H}$,
$L_{45} = 0.3655 \, \mathrm{H}$, $L_{46} = 1.0245 \, \mathrm{H}$,
$L_{36} = 0.3240 \, \mathrm{H}$, $L_{67} = 0.4735 \, \mathrm{H}$,
$L_{37} = 0.1875 \, \mathrm{H}$, $L_{47} = 0.74 \, \mathrm{H}$. 

\subsubsection{No Guarantee on Synchronization}
$L_{12} = 3.336 \, \mathrm{H}$, $L_{15} = 2.684 \, \mathrm{H}$,
$L_{26} = 1.108 \, \mathrm{H}$, $L_{56} = 4.23 \, \mathrm{H}$,
$L_{45} = 1.462 \, \mathrm{H}$, $L_{46} = 4.098 \, \mathrm{H}$,
$L_{36} = 1.296 \, \mathrm{H}$, $L_{67} = 1.894 \, \mathrm{H}$,
$L_{37} = 0.75 \, \mathrm{H}$, $L_{47} = 2.96 \, \mathrm{H}$.
\end{appendix}
\bibliographystyle{ieeetr}
\addcontentsline{toc}{section}{\refname}\bibliography{references}
\end{document}